\documentclass[a4paper,12pt]{article}

\usepackage[english]{babel}
\usepackage[tbtags]{amsmath}
\usepackage{amssymb}
\usepackage{amsfonts}
\usepackage{amsthm}
\usepackage{calrsfs}
\usepackage{array}
 \usepackage{bbm}
\usepackage{stmaryrd}
\usepackage{upgreek}
\usepackage[svgnames]{xcolor}
\usepackage{hyperref}
\usepackage{mathtools}
\def\div{\, \mbox{div}\,  }

\def\pp{ p-\varepsilon}

\def\p{\phi }

\numberwithin{equation}{section}

\def\iint {\int_B}

\def\grad{\partial_y}

\def \I {B(x_0,T(x_0)-t)}

\def\Box{\hfill\rule{2.5mm}{2.5mm}}
\def\t{{\mathrm{d}}{\tau}}

\def\y{{\mathrm{d}}y}
\def\e{\varepsilon}
\def\no{\nonumber}

\def\grad{\nabla}

\def \er {\mathbb R}

\def\R{{\mathbb {R}}}

\def\cprime{$'$}

\newcommand{\ds}{\displaystyle}
\newcommand{\dv}{\mathop{div}}

\theoremstyle{remark}
\newtheorem{nb}{Remark}[section]

\makeatletter
\def\blfootnote{\xdef\@thefnmark{}\@footnotetext}
\makeatother

\title{\bf The blow-up rate for    a non-scaling invariant   semilinear wave equations  in higher dimensions}

\date{\today}

\allowdisplaybreaks

\def\Box{\hfill\rule{2.5mm}{2.5mm}}
\def\t{{\mathrm{d}}\tau}

\def\y{{\mathrm{d}}y}
\def\e{\varepsilon}
\def\no{\nonumber}

\def\no{\nonumber}

\def \no { \noindent}
\def \er {\mathbb R}

\def\t0{{t}}

\def \p {\rho_{\varepsilon}}

\def\grad{\nabla}
\def\div{\, \mbox{div}\,  }

\def\MOD#1{{|\kern -.16em |\kern -.16em | #1 | \kern -.16em |\kern
 -.16em |}}
\def \epsilon {\varepsilon}

\def\ds{\displaystyle}
\newtheorem{theo}{\bf THEOREM}
\newtheorem{coro}[theo]{\bf COROLLARY}
\newtheorem{lem}{\bf LEMMA}[section]
\newtheorem{pro}[lem]{\bf PROPOSITION}
\newtheorem{cor}[lem]{\bf COROLLARY}

\newtheorem{rem}[lem]{\bf REMARK}

\numberwithin{equation}{section}

\def\B1{B_{1/2}}
\def\Box{\hfill\rule{2.5mm}{2.5mm}}

\def\iint{\int_{B}}

\def\N{{\mathbb {N}}}
\def\R{{\mathbb {R}}}

\def\build#1_#2^#3{\mathrel{
\mathop{\kern 0pt#1}\limits_{#2}^{#3}}}

\def\dv{\mathop{\rm div}}

\def\h1{\mathop{\rm H^1_{\rm loc,\rm u}}}

\def\l2{\mathop{\rm L^2_{\rm loc,\rm u}}}

\def\t0{{t}}

\def \pp {\frac{\rho_{\e}}{\sqrt{1-|y|^2}}}
\def\Box{\hfill\rule{2.5mm}{2.5mm}}

\author{Mohamed  Ali  Hamza\\
{\it \small 
Imam Abdulrahman Bin Faisal University,
 Dammam, 34212, Saudi Arabia}\\
Hatem Zaag\\
{\it \small Universit\'e Sorbonne Paris Nord},\\
{\it \small LAGA, CNRS (UMR 7539), F-93430, Villetaneuse, France}
}

\begin{document}
\bibliographystyle{siam}

\title{\textbf{A better bound on blow-up rate for the superconformal semilinear wave equation}}


\maketitle

%


\begin{abstract} We consider
the semilinear wave equation in higher dimensions with superconformal
power nonlinearity.
The purpose of this paper is to give a new upper bound on the blow-up
rate in some space-time integral, showing a $|\log(T-t)|^q$
improvement in comparison with previous results obtained in
\cite{HZdcds13,KSVsurc12}. 

%
\end{abstract}

\vspace{0.5cm}


 \vspace{0.4cm}

\section{\bf Introduction}

\subsection{Motivation of the problem}

This paper is devoted to the study of blow-up solutions for the
following semilinear  wave equation:
\begin{equation}\label{gen}
\left\{
\begin{array}{l}
\partial_t^2 u =\Delta u+|u|^{p-1}u,\\
\\
(u(x,0),\partial_tu(x,0))=(u_0(x),u_1(x))\in  H^{1}_{loc}(\er^N)\times L^{2}_{loc}(\er^N),\\
\end{array}
\right.
\end{equation}
in space dimension $N\ge 2$, where $u(t):x\in \er^N \rightarrow
u(x,t)\in \er$,
\begin{equation}\label{surc}
p_c<p<p_S,
\end{equation}
 $p_c\equiv 1+\frac4{N-1}$ is the conformal critical exponent and  $p_S\equiv 1+\frac4{N-2}$ if $N\ge 3$
is the Sobolev critical exponent.

\bigskip

 \no The Cauchy problem of equation ({\ref{gen}}) is solved in
$H^{1}_{loc}\times L^{2}_{loc}$. This follows from the finite speed
of propagation and the wellposedness in $H^{1} \times L^{2}$,
valid whenever $ 1< p< p_S$. 
Combining the blow-up solutions for the associated ordinary differential equation of (\ref{gen})   and 
the finite speed of
propagation, we conclude that there exists a
 blow-up solution $u(t)$ of (\ref{gen})  which depends non-trivially on the space variable. 
Numerical simulations of blow-up are given by Bizo\'n {\it et al.} in \cite{BBMWnonl10,BCSjmp11, BCTnonl04}.

\medskip

If $u$ is an arbitrary blow-up solution of   \eqref{gen}, we define
a 1-Lipschitz curve $\Gamma=\{(x,T(x))\}$
such that
the domain of definition of $u$ is written as
\begin{equation}\label{defdu}
D=\{(x,t)\;|\; t< T(x)\}.
\end{equation}
The justification of this fact can be found in Alinhac
\cite{Apndeta95}, where such a set is referred to as the \textit{maximal influence domain}. 
 In fact, for any $x_0\in \R^N$, we derive a solution
which is defined in $K_{x_0,t_0}$, the backward light cone
defined by
\[
K_{x_0,t_0}=\{(\xi,\tau)\in \R^N\times[0,t_0)\;|\;|\xi-x_0|<t_0-\tau\}.
\]
Denoting by $T(x_0)$ the supremum of such a $t_0$, we fall into one of
these two situations:\\
- either for some $x_1\in \R^N$, $T(x_1)=+\infty$, and in that case,
the solution if defined for all $x\in \R^N$ and $t\ge 0$. In other
words, $T(x_0)=+\infty$, for any $x_0\in \R^N$ ;\\
- or, for all $x_1\in \R^N$, $T(x_1)<+\infty$, and this is the case we
are interested in. In particular,  the domain of definition of $u$
is the union of backward light cones of the type $K_{x_0,t_0}$, and
such a union is clearly of the form \eqref{defdu} with $x\mapsto T(x)$
is $1-$Lipschitz. See also Sasaki \cite{S18ade} for a related definition of the
domain of definition in some similar setting.

\medskip

The time
$\bar T=\inf_{x\in {\R^N}}T(x)$ and $\Gamma$ are called the blow-up
time and the blow-up graph of $u$.
A point $x_0$ is non characteristic
if  there are
\begin{equation}\label{nonchar}
\delta_0\in(0,1)\mbox{ and }t_0<T(x_0)\mbox{ such that }
u\;\;\mbox{is defined on }{\mathcal C}_{x_0, T(x_0), \delta_0}\cap \{t\ge t_0\},
\end{equation}
where ${\cal C}_{\bar x, \bar t, \bar \delta}=\{(x,t)\;|\; t< \bar
t-\bar \delta|x-\bar x|\}$. If not, $x_0$ is said to be characteristic.

\bigskip

When the exponent $p$ is conformal or subconformal, i.e. when 
\begin{equation}\label{sousc1}
1<p\le p_c,
\end{equation}
Merle and Zaag proved in \cite{MZajm03,MZma05,MZimrn05} that the blow-up rate of equation \eqref{gen} is $C(T-t)^{-\frac 2{p-1}}$, the solution of the associated ODE $u"=u^p$. 
 More precisely, 
they proved  that if $u$  is a solution of \eqref{gen}
with blow-up graph $\Gamma : \{ x\mapsto T(x)\}$ and $x_{0}$ is a
non-characteristic point, then,  for all $t\in
[\frac{3T(x_{0})}{4},T(x_{0})]$,
 \begin{eqnarray}\label{mzmz}
  0 < \varepsilon_{0}(p)\leq (T(x_{0})-t)^{\frac{2}{
   p-1}}\frac{\|u(t)\|_{L^{2}(\I )}}
{(T(x_{0})-t)^{\frac{N}2}}\\
  +(T(x_{0})-t)^{\frac{2}{ p-1}+1}\Big(\frac{\|\partial_{t} u(t)\|_{L^{2}(\I)}}
  {(T(x_{0})-t)^{\frac{N}2}}
  + \frac{\|\partial_x u(t)\|_{L^{2}(\I )}}{{(T(x_{0})-t)^{\frac{N}2}}}\Big)\leq K_0,\nonumber
  \end{eqnarray}
 where the constant $K_0$ depends only on  $p$ and on an upper bound on
 $T(x_{0})$, ${1}/{T(x_{0})}$, $\delta_{0}(x_{0})$, together with the
 norm of initial data
 in $H^{1}_{loc,u}(\R^N)\times L^{2}_{loc,u}(\R^N)$. 
Their method relies on the existence of a Lyapunov functional in self-similar variables 
 which  restricts their results to subconformal and conformal range  (i.e. $1<p\le \frac{N+3}{N-1})$. In fact,
 when $p>\frac{N+3}{N-1}$, there is no such a Lyapunov functional, and the method breaks-down.

\medskip

In \cite{HZjhde12, HZnonl12}, using a highly non-trivial perturbative method, we could obtain the blow-up rate for the
 Klein-Gordon equation  and more generally, for equation 
\begin{equation}\label{NLWP}
\partial_t^2 u =\Delta u+|u|^{p-1}u+f(u)+g(\partial_t u ),\,\,\,(x,t)\in \R^N \times [0,T),
 \end{equation}
under the assumptions  $|f(u)|\leq M(1+|u|^q)$ and $|g(v)|\leq M(1+|v|)$,   for some $M > 0$ and $q<p\le  \frac{N+3}{N-1}$. 
In fact, we proved a similar result to $\eqref{mzmz}$,
valid in the subconformal  and conformal case.
  Let us also mention that in  \cite{H1, omar1, omar2}, the
results obtained in  \cite{HZjhde12,HZnonl12}   were extended by Hamza and Saidi to the  
strongly perturbed equation  (\ref{NLWP})  with  $|f(u)|\leq
M(1+|u|^p\log^{-a}(2+u^2))$,  for some $a > 1$, though keeping the
same condition in $g$,  valid in the subconformal  and conformal
case. 

\medskip

We would like to mention the result of Donninger and Sch\"orkhuber in
\cite{DSdpde12} in 3 space dimensions in the subcritical range
(i.e. $1<p\le 3)$, who proved  the stability of the ODE solution 
$u(t)=\kappa_0(p)(T-t)^{-\frac 2{p-1}}$ among all radial solutions,
with respect to small perturbations in initial data in the energy
space. 
Later, they extended their results to the superconformal range
(i.e. $p>3)$, though, in a topology stronger than the energy
space. They did it first in the radial case in \cite{DStams14}, then
without any symmetry assumption in \cite{DScmp16}. 
 Their approach is based in particular on a good understanding of the
 spectral properties of the linearized operator in self-similar
 variables, which happens to be non-self-adjoint.
 Later, by  establishing 
 suitable Strichartz estimates for the critical wave equation
 (i.e. $p=3$)     in similarity variables,  
 Donninger proved in \cite{DDuke}
 the asymptotic stability of the ODE  for the  radial solution  
with respect to small perturbations in initial data, in the energy space.

\bigskip

Let us focus now on the superconformal and Sobolev subcritical range
\eqref{surc}, i.e. when  
\begin{equation}\label{psurconforme}
p_c<p<p_S.
\end{equation}
To the best of our knowledge, the first available result is due to
Killip, Stovall and Vi\c san in \cite{KSVsurc12}, who give an upper
bound on the blow-up rate for the Klein-Gordon equation, among many
other interesting results. More precisely, if $u(x,t)$ is a blow-up solution of equation \eqref{NLWP} with $(f,g)\equiv (-u,0)$ and $\{x\mapsto T(x)\}$ is its blow-up graph, then  for all  $x_0\in \er^N$, there exists $K_1>0$ such that, for all $t\in [0,T(x_0))$,
\begin{equation}\label{v1}
(T(x_0)-t)^{-\frac{(p-1)N}{p+3}}\int_{B(x_0,\frac{T(x_0)-t}2)} \! u^2(x,t)\mathrm{d}x\;\; \le\;\;
K_1,
\end{equation}
and for all   $t\in (0,T(x_0)]$,
\begin{equation}\label{v2}
\int_{T(x_0)-t}^{T(x_0)-\frac{t}2}\int_{B(x_0,\frac{T(x_0)-\tau }2)}
\!\Big(|\grad u(x,\tau )|^{2}+|\partial_tu(x,\tau
)|^2\Big)\mathrm{d}x\mathrm{d}\tau\;\;\le\;\;K_1.
\end{equation}
Moreover, if $x_0$ is non-characteristic,
then estimates \eqref{v1} and \eqref{v2} hold with the ball $B\left(x_0,\frac{T(x_0)-\tau}2\right)$ replaced by the ball $B(x_0,T(x_0)-\tau )$.
Subsequently, in   \cite{HZdcds13}, thanks to a different method based 
on the use of self-similar variables, we refined estimate
\eqref{v1} for any point  $x_0\in \R^N$ by showing that
\begin{equation}\label{vv1}
(T(x_0)-t)^{-\frac{(p-1)N}{p+3}}\int_{B(x_0,{T(x_0)-t})} \! u^2(x,t)\mathrm{d}x \rightarrow 0 ,\ \ \textrm{ as}\ \ \ t\rightarrow T(x_0).
\end{equation}
 Furthermore, if $x_0$ is a non-characteristic point, we refined
 \eqref{v2} by proving that
\begin{equation}\label{vv2}
\int_{T(x_0)-t}^{T(x_0)-\frac{t}2}\int_{B(x_0,T(x_0)-\tau )}
\!\Big(|\grad u(x,\tau )|^{2}+|\partial_tu(x,\tau
)|^2\Big)\mathrm{d}x\mathrm{d}\tau\rightarrow 0 ,\ \ \textrm{ as}\ \ \ t\rightarrow T(x_0).
\end{equation}

 \bigskip

Let us mention that both our approach in  \cite{HZdcds13} and
that of \cite{KSVsurc12} rely on the existence of a Lyapunov
functional. While the authors of \cite{KSVsurc12} find it in the
$u(x,t)$ setting, we find a different one, in the $w(y,s)$
setting. The key difference lays in the dissipation of the Lyapunov
functional, which involves an unusual term in our approach, namely
$\int |w|^{p+1}$, unlike the functional of of \cite{KSVsurc12}. Since
the dissipation has to go to zero in a certain sense, this enables us
to show that some $L^q$ norms of $w$ have to be negligible with respect
to the rate predicted by those authors, and not just bounded by that rate. More precisely, our method  in  \cite{HZdcds13} shows that some quantities go to zero (see \eqref{vv1} and \eqref{vv2}), while \cite{KSVsurc12} only shows their boundedness (see \eqref{v1} and \eqref{v2}).

\bigskip

In this paper,  our purpose is to show a new upper bound on the blow-up rate which is of course an improvement of the previous results.
More precisely, we    state  the following:

\begin{theo}\label{t0}\textbf{(Growth estimate near the blow-up surface for equation ({\ref{gen}}))}.\\
Consider  $u$  a solution of equation ({\ref{gen}}) with blow-up graph
$\{x\mapsto T(x)\}$, and $q>0$. Then, if   $x_0\in \er^N$   is   a non-characteristic point,   there exists $K>0$, such that   for all  $t\in
[0,T(x_0))$,   we have
\begin{equation}\label{t4}
|\log (T(x_0)-t)|^q\int_t^{\frac{t+T(x_0)}2}\int_{B(x_0,{T(x_0)-\tau })}
\!\Big(|\grad u(x,\tau )|^{2}+|\partial_tu(x,\tau
)|^2\Big)\mathrm{d}x\mathrm{d}\tau\le K.
\end{equation}
%
%
%
%
Moreover, we have
\begin{eqnarray}\label{limener}
&&|\log (T(x_0)-t)|^q\frac{T(x_0)-t}2\int_{B(x_0,{T(x_0)-t
})} \!\Big( |\grad u(x,t )|^{2}- \big(\frac{x-x_0}{T(x_0)-t}\cdot \grad u(x,t )\big)^2\nonumber\\
&&+|\partial_tu(x,t )|^2-\frac1{p+1} \!|
u(x,t)|^{p+1}\Big)\mathrm{d}x\le K.
\end{eqnarray}
\end{theo}
Thanks to \eqref{t4}, and \eqref{limener}, we easily refine estimate \eqref{vv1}. More precisely, we show the following:
\begin{coro}
Under the hypotheses of Theorem \ref{t0},  we have
\begin{equation}\label{vv1bb}
  |\log (T(x_0)-t)|^q(T(x_0)-t)^{-\frac{(p-1)N}{p+3}}\int_{B(x_0,{T(x_0)-t})} \! u^2(x,t)\mathrm{d}x \le \tilde{K}.
\end{equation}
\end{coro}

\begin{nb}
The constant $K$, and the rate of convergence to $0$ of the
different quantities in the previous theorem  and in the whole
paper, depend only on $N$, $p$, $q$ and some upper bound on $T(x_0)$,
$1/{T(x_0)}$, and the initial data $(u_0,u_1)$ in $
H^{1}(B(x_0,2T(x_0)))\times L^{2}(B(x_0,2T(x_0)))$, together with
$\delta_0(x_0)$ if $x_0$ is non-characteristic point.
\end{nb}

\begin{nb}
  Let us remark that we have the following   lower bound which
follows  from standard techniques (scaling arguments, the
wellposedness in $H^1(\er^N )\times L^2(\er^N)$, the finite speed of
propagation and the fact that $x_0$ is a non-characteristic point):
there exist
 $\varepsilon_0>0$ such that
\begin{eqnarray*}
&&0<\varepsilon_0\le (T(x_0)-t)^{\frac{2}{p-1}}\frac{\|u(t)\|_{L^2(B(x_0,{T(x_0)-t}))}}{ (T(x_0)-t)^{\frac{N}{2}}}\nonumber\\
&&+ (T(x_0)-t)^{\frac{2}{p-1}+1}\Big
(\frac{\|\partial_tu(t)\|_{L^2(B(x_0,{T(x_0)-t}))}}{
(T(x_0)-t)^{\frac{N}{2}}}+
 \frac{\|\grad u(t)\|_{L^2(B(x_0,{T(x_0)-t}))}}{ (T(x_0)-t)^{\frac{N}{2}}}\Big
 ).
\end{eqnarray*}
For details, see the argument on page 1149 of \cite{MZajm03}. Indeed,
although that paper was written in the subconformal rage, the argument
holds for all Sobolev subcritical exponent.
\end{nb}
\begin{nb}
Since we crucially need a covering technique in our argument (in
particular, in the construction of some family of Lyapunov
functionals), our method breaks down too in the case of a characteristic point and we are unable to obtain similar estimates.
\end{nb}
\begin{nb}
 Up to some time-dependent factor, the expression in (\ref{limener})
is equal to the main terms of the energy  in similarity variables, and this is an improvement with respect to 
previous results 
obtained in \cite{HZdcds13,KSVsurc12} in the superconformal range. 
However, even with this improvement, we think that our estimates are still not optimal. 
\end{nb}

\begin{nb}
Let us denote that,
following the analysis of  Hamza in Saidi in \cite{omar1}, our result holds    for 
more general   equations of the type
\begin{equation}\label{ABCD}
\partial_t^2 u =\Delta u+{f}(u)+g(x,t,u,\partial_t u,\grad u ),\,\,\,(x,t)\in \R^N \times [0,T),
 \end{equation}
where $ {f}:\er \to \er$ and  $ {g}:\er^{2N+3}
\to \er$ are ${\cal C}^1$ functions which satisfy the following conditions:
\begin{equation}\label{H1}
|g(x,t,u,v,w)|\leq \tilde C(1+|u|^{\frac{p+1}2}+|v|+|w|),   
 \end{equation}
and 
\begin{equation}\label{H2}
\Big|\int_0^u{f}(y)
{\mathrm{d}}y
-\frac{u{f} (u)}{p+1}\Big|\le \tilde C\Big(1+  |u|^{p+1} \log^{-a}(2+u^2)\Big),
\end{equation}
for some $\tilde C > 0$, $p>1$,  and $a>1$. 
Indeed,   under these assumptions,  
we can construct  a suitable Lyapunov
functional for equation \eqref{ABCD} under hypotheses  \eqref{H1} and \eqref{H2}. Therefore, we  can show  similar results to \eqref{t4} and \eqref{limener}.  
However,
In order to keep our analysis clear, we  restrict our analysis to   equation \eqref{gen}.
\end{nb}

\subsection{Strategy of the proof}

In this subsection, we would like to briefly present our method to prove Theorem \ref{t0}. It relies on the estimates in similarity variables
introduced in \cite{AMimrn01} and used in \cite{MZajm03,MZimrn05,MZma05}. More precisely, given $(x_0,T_0)$
such that $0< T_0\le T(x_0)$, we introduce  the following
self-similar change of variables:
\begin{equation}\label{scaling}
y=\frac{x-x_0}{T_0-t},\qquad s=-\log (T_0-t),\qquad
u(x,t)=(T_0-t) ^{-\frac{2}{p-1}}w_{x_0,T_0}(y,s).
\end{equation}
If $T_0=T(x_0)$, we may write $w_{x_0}$ instead of $w_{x_0,T(x_0)}$,
for simplicity.
From equation (\ref{gen}), the  function $w_{x_0,T_0}$  (we write $w$ for
simplicity) satisfies the following equation for all $y\in B\equiv
B(0,1)$ and $s\ge -\log T_0$:
\begin{eqnarray}\label{C}
\partial_{s}^2w&=&\div( \grad w-(y\cdot \grad w)y)-2\alpha
 y\cdot\grad w
-\frac{2p+2}{(p-1)^2}w+|w|^{p-1}w\\
&&-\frac{p+3}{p-1}\partial_s w-2y\cdot \grad \partial_sw,\
  \forall y\in B\ {\textrm{ and}} \ s\ge -\log T_0,\nonumber
\end{eqnarray}
where
\begin{equation}\label{alpha}
\alpha=\frac2{p-1}-\frac{N-1}{2}=\frac2{p-1}-\frac2{p_c-1}<0,
\end{equation}
%
In particular, the behavior of $u$ as $t \rightarrow T_0$ and that of
$w$ as $s \rightarrow +\infty$ are obviously related. 

\medskip

We shall now    notice  that 
in the conformal case where $p=p_c$, Merle and Zaag \cite{MZimrn05}
proved that
\begin{equation}\label{f} E_0(w(s))=\displaystyle\int_{B}\Big (
\frac{1}{2}(\partial_sw)^2 +\frac{1}{2}|\grad
w|^2-\frac{1}{2}(y\cdot \grad w)^2+\frac{p+1}{(p-1)^2}w^2
-\frac{|w|^{p+1}}{p+1}\Big ) {\mathrm{d}}y,\quad
\end{equation}
 is a Lyapunov functional for equation (\ref{C}).
Furthermore,  when $p>p_c$, we showed that 
 the  following energy functional
\begin{eqnarray}
\label{00E3}
 F_0(w(s),s)&=& {E(w(s))e^{2\alpha s}},
\end{eqnarray}
where  $\alpha$ is defined in \eqref{alpha}, and  
\begin{eqnarray}\label{E1}
E(w(s))&=&\displaystyle {E_0(w(s))+\alpha\displaystyle\int_{B}w\partial_s w {\mathrm{d}}y- \frac{\alpha N}2 \displaystyle\int_{B}w^2
 {\mathrm{d}}y},
\end{eqnarray}
is also a Lyapunov functional for equation (\ref{C}).  This  allows us just
to get
a rough exponential space-time estimate of the solution $u$   of \eqref{gen} obtained in 
\cite{HZdcds13} in similarity variables. More precisely, we
established the following results:

\bigskip

\noindent {\it{\bf { (Exponential bound on space-time integrals  of solutions to
(\ref{gen}) in similarity variables.}} Let  $u$ be  a solution of (\ref{gen}) with
blow-up graph $\Gamma :\{x\mapsto T(x)\}$. Then,
 for all  $x_0\in \er^N$  and  $T_0\le
T(x_0)$, we have
 for all $s\ge s_0=-\log T_0$, 
\begin{eqnarray}\label{cor3}
e^{\frac{8\alpha s}{p+3}}\int_{ B}\!|w(y,s)|^2{\mathrm{d}}y \le  K_2,
\end{eqnarray}

\begin{equation}\label{F0}
\int_{s_0}^{\infty}e^{2\alpha s}
\!\!\int_{B}\!|w(y,s)|^{p+1}{\mathrm{d}}y{\mathrm{d}}s \le
 K_2.
\end{equation}
If in addition $x_0$ is non-characteristic (with a slope $\delta_0\in (0,1))$, then
\begin{equation}\label{FF1}
e^{2\alpha s}\int_{s}^{s+1} \!\int_{B}\!\big((\partial_s
w(y,\tau))^2
+|\grad w(y,\tau)|^2\big)
{\mathrm{d}}y{\mathrm{d}}\tau\ \ \ \rightarrow 0
\textrm{ as}\ \ \ s\rightarrow +\infty,
\end{equation}
 where the constant $K_2$ depends only on $N, p, \delta_{0}(x_{0})$, 
 $T(x_{0})$ and 
$\|(u_0,u_1)\|_{H^{1}\times
  L^{2}(B(x_0,\frac{T(x_0)}{\delta_0(x_0)}) )}$.}

\bigskip

\noindent Furthermore, for all $s'\ge s\ge -\log T_0=s_0$, the functional
$F_0(w(s),s)$ defined in (\ref{00E3}) satisfies
\begin{eqnarray}\label{LE0}
F_0(w(s'),s')-F_0(w(s),s) &=&
- \int_{s}^{s'} e^{2\alpha \tau}\!\int_{\partial B}\!\Big(\partial_s w+\alpha w\Big)^2{\mathrm{d}}\sigma {\mathrm{d}}\tau\nonumber\\
&&+ \frac{\alpha (p-1)}{p+1}\int_{s}^{s'} e^{2\alpha
\tau }\int_{B}|w|^{p+1}{\mathrm{d}}y{\mathrm{d}}\tau.
\end{eqnarray}
Moreover, for all $s\ge s_0$, we have  $F_0(w(s),s)\ge0$, and 
\begin{equation}\label{nesF0}
F_0(w(s),s)\rightarrow 0
\textrm{ as }s\rightarrow +\infty.
\end{equation}
.


In fact, we believe that the good blow-up rate in the superconformal
case $p>p_c$ (still in the Sobolev subcritical range) is given by the
associated ODE solution $u''=u^p$, as it is already the case for $p\le
p_c$. In other words, this can be expressed as the boundedness of $w(y,s)$, the
similarity variables' version defined in \eqref{scaling}. A key estimate towards
such an estimate would require bounding  the functional
 $E(w(s))=e^{-2\alpha s}F_0(w(s),s)$. Unfortunately, we have no
 satisfactory answer to this prediction.  Nonetheless, we make in this
 paper a step towards bounding that quantity, by improving estimate
 \eqref{nesF0}.  

 \medskip
 
In order to do so, the first important step consists in studying two
functionals associated to equation \eqref{C}, the similarity
variables' version of equation \eqref{gen}. This will be done in the
non-characteristic case.

\medskip

The first functional (called $N_{\frac12+\e}(w(s))$ below, where
$\e>0$ is arbitrary) arises from a Pohozaev identity obtained through 
the multiplication of equation \eqref{C} by $y.\grad
w\sqrt{1-|y|^2}\p$ where
\begin{equation}\label{A50}
\p(y)=(1-|y|^2)^{\varepsilon}.
\end{equation}
The second functional is obtained by multiplying
equation \eqref{C} by $w(1-|y|^2)^{-\frac12}\p$. By combining these functionals, we easily derive a new estimate controling 
the time average of the  $L^{p+1}(B \times [s,s+1])$ norm of $w$ with the singular weight $\frac{\p}{\sqrt{1-|y|^2}}.$
This estimate is given in Proposition \ref{P01m} below, which is one of the novelties of this paper.

\bigskip

%
%
%
%

Using the estimates mentioned above, 
together with  the  boundedness of the solution $w$
in the space-time norm $L^{p+1}(B \times [-\log T,+\infty))$, we 
get a better estimate of the $\dot H^1(B \times [-\log T,+\infty))$ norm  
of   the solution of equation \eqref{C}.  More precisely, we obtain estimates
\eqref{K2} below, and thanks to  a  covering argument,  we deduce that
\begin{equation*}
\int_{\tilde s_0}^{\infty}\tau^{\nu_0-1}e^{2\alpha \tau}
 \int_{B}\Big((\partial_s w)^2 + |\grad w|^2\Big) \y
{\mathrm{d}}\tau  \le C,
\end{equation*}
for some $\nu_0\in (0,1)$.
 Hence, by using the above estimate together with\eqref{LE0}, we derive that the functional $s^{\nu_0}F_0(w(s),s)$ is also a  Lyapunov functional for
equation \eqref{C} which allows us to write
\begin{equation}\label{nesF0bis}
s^{\nu_0}F_0(w(s),s)\rightarrow 0
\textrm{ as } s\rightarrow +\infty.
\end{equation}
With this new
Lyapunov functional  $s^{\nu_0}F_0(w(s),s)$ at hand, the adaptation of the interpolation
strategy from our previous papers works straightforwardly.
Consequently, we show that
\begin{equation*}
s^{\nu_0}e^{2\alpha s}\int_{s}^{s+1} \!\int_{B}\!\big((\partial_s
w(y,\tau))^2
+|\grad w(y,\tau)|^2+|w(y,s)|^{p+1}\big)
{\mathrm{d}}y{\mathrm{d}}\tau\ \ \ \rightarrow 0
\textrm{ as}\ \ \ s\rightarrow +\infty.
\end{equation*}
 Furthemore, by an {\bf {induction argument}} we prove that  for any   $k\in \N,$  the functional
 $s^{k\nu_0}F_0(w(s),s)$ is also a Lyaponov functional, and that
\begin{equation}\label{nesFknu}
s^{k\nu_0}F_0(w(s),s)\rightarrow 0
\textrm{ as }s\rightarrow +\infty.
\end{equation}
Similarly, relying on the family of new Lyapunov functionals
$\Big(s^{k\nu_0}F_0(w(s),s)\Big)_{k\in \N}$ and following the same
strategy as in the particular case $k=1$, we show that for all $k\in \N$,
\begin{equation}\label{nesFknu000}  
s^{k\nu_0}e^{2\alpha s}\int_{s}^{s+1} \!\int_{B}\!\big((\partial_s
w(y,\tau))^2
+|\grad w(y,\tau)|^2+|w(y,s)|^{p+1}\big)
{\mathrm{d}}y{\mathrm{d}}\tau\ \ \ \rightarrow 0
\textrm{ as}\ \ \ s\rightarrow +\infty.
\end{equation}
 Obviously, estimate \eqref{nesFknu000} implies that for any $k\in \N$,
 \begin{equation}\label{nesFknu000b}  
\int_{s_0}^{\infty}\tau^ke^{2\alpha \tau} \!\int_{B}\!\big((\partial_s
w(y,\tau))^2
+|\grad w(y,\tau)|^2+|w(y,s)|^{p+1}\big)
{\mathrm{d}}y{\mathrm{d}}\tau\ \ \ \rightarrow 0
\textrm{ as}\ \ \ s\rightarrow +\infty.
\end{equation}
The implication in the other way is even more obvious.


\bigskip

As we mentioned before, our first novelty  in this paper lays in the
following proposition:  
\begin{pro}\label{P01m}
For all $\e>0$,  for all $s\geq -\log(T_{0})+2$ we have 
\begin{equation}\label{w01}
\int_{s}^{s+1}\int_{B}|w|^{p+1}\frac{\p}{\sqrt{1-|y|^2}}\y 
{\mathrm{d}}\tau \le C\int_{s-2}^{s+3}\!\int_{B}\Big(|\nabla  w|^2+(\partial_sw)^2+w^2+|w|^{p+1} \Big) \p {\mathrm{d}}y{\mathrm{d}}\tau. 
\end{equation}
%
\end{pro}

In the sequel, we will use the following notations for the radial
and the angular derivatives:
\begin{equation}\label{wr}
 \nabla_r w=\frac{y\cdot\nabla w}{|y|^2} y\quad \ {\textrm{and}}\quad  \nabla_{\theta}w=
\nabla w-
\frac{y\cdot\nabla w}{|y|^2} y.
\end{equation}
Clearly, it holds that $\nabla w=\nabla_rw+\nabla_{\theta}w$ and from
orthogonality between $\nabla_rw$ and $\nabla_{\theta}w$, it follows that
\begin{equation}\label{wr1}
|y|^2|\grad w|^2-(y\cdot \grad w)^2=|y|^2|\nabla_{\theta}w|^2,
\end{equation}
\begin{equation}\label{wr2}
(y\cdot \grad w)^2 =|y|^2|\nabla_{r}w|^2,
\qquad \textrm{
and}\qquad 
|\grad w|^2=|\nabla_{\theta}w|^2+|\nabla_{r}w|^2.
\end{equation}
%
%
%
%
Note also that equation (\ref{C}) will be studied in the Hilbert space $\cal H$
\begin{equation}\label{HH}
{\cal H}=\Big \{(w_1,w_2), |
\displaystyle\int_{B}\Big ( w_2^2 +|\grad w_1|^2-(y\cdot \grad w_1)^2)+w_1^2\Big) {\mathrm{d}}y<+\infty \Big \}.
\end{equation}
%
%
Throughout this paper, $C$ denotes a generic positive constant
 depending only on $p$  and $N,$  which may vary from line to line.
 In addition, we  will use $K_1,K_2,K_3...$  as  positive constants
 depending only on $p,N, \delta_0(x_0)$, the slope definind a
 non-characteristic point (see \eqref{nonchar}), as well as initial data, which
 may also vary from line to line.  We will also write $f(s)\sim g(s)$
 to indicate that
 $\displaystyle{\lim_{s\to \infty}\frac{f(s)}{g(s)}=1}.$ 

\bigskip

The paper is organized as follows:  In Section \ref{sec2}, we prove
Proposition \ref{P01m}.   In Section \ref{sec3}, we give an
improvement of the bound to the $\dot{H}^{1}(\R^N\times
[-\log(s_0),\infty))$ norm of $w$, based on the results.  This allows
us to construct a family of  Lyapunov functionals for equation
\eqref{C}. Finally, in Section \ref{sec4} and subsequently, we
establish Theorem \ref{t0}. \\

\par 

{\bf{Acknowledgement.}} 
Hatem Zaag wishes to thank Pierre Raphaël and the ”SWAT” ERC project for their support.

\section{A bound to the time average of the $L^{p+1}$ norm of $w$ with singular weight}\label{sec2}

This section is devoted 
to the statement and the proof of a general version of Proposition
\ref{P01m}, uniform for $x$ near $x_0$ (see Proposition \ref{P01m}' below). This
section is divided into three subsections:
\begin{itemize}
\item In the first one, we give some classical energy estimates
following from the multiplication of equation \eqref{C} by
$w(1-|y|^2)^{\e}$ and $\partial_sw(1-|y|^2)^{\e}$, then integration,
for any $\e\in (0,1)$. 
\item The second subsection is devoted to  a new energy estimate
  following from a Pohozaev multiplier. More precisely, we multiply
  equation \eqref{C} by $y\cdot\grad w(1-|y|^2)^{\e}$ and
  $w(1-|y|^2)^{-\frac12+\e}$,  for any $\e\in (0,1)$, then integrate, in order to get
  a new identity.
\item  By combining the above energy estimates,  we conclude the proof of  Proposition \ref{P01m}'.
\end{itemize}

Let us first recall the following uniform version of the exponential
bounds on space-time integral of the solution $u$   of \eqref{gen} near any
non characteristic point  obtained in 
\cite{HZdcds13} and written in similarity variables \eqref{scaling}:

\medskip

\noindent {\it{\bf { (Uniform exponential bounds on space-time
      integrals in similarity variables).}}}
{\it
\noindent 
Let $u $   a solution of ({\ref{gen}}) with
blow-up graph $\Gamma:\{x\mapsto T(x)\}$ and  $x_0$  a non
characteristic point. Then,   for all $T_0 \in (0,T(x_{0})]$,  for all $s\geq -\log T_0$ and $x\in \er^N$ where $|x-x_0|\le \frac{T_0-t}{\delta_0(x_0)}$, 
we have
\begin{equation}\label{feb19}\int_s^{s+1}\int_B \big( |\grad w|^2 +( \partial_{s}w)^2\big)\y d\tau \leq K_3 e^{-2\alpha s},
\end{equation}
\begin{equation}\label{feb191}
\int_{s_0}^{\infty}e^{2\alpha \tau}\int_B | w(y,\tau)|^{p+1}\y d\tau \leq K_3,
\end{equation}

\begin{eqnarray}\label{cor3t}
e^{\frac{8\alpha s}{p+3}}\int_{ B}\!|w(y,s)|^2{\mathrm{d}}y \le K_3,
\end{eqnarray}
where $w=w_{x,T^*(x)}$ is defined in \eqref{scaling}, with
 \begin{equation}\label{18dec1}
T^*(x)=T_0-\delta_0(x_0)|x-x_0|,
\end{equation}
 and $\delta_{0}(x_{0})$ defined in \eqref{nonchar}.  Note that  $K_3$ depends on $ p,  N,  \delta_{0}(x_{0})$, 
 $T(x_{0})$,  and
$\|(u_0,u_1)\|_{H^{1}\times
L^{2}(B(x_0,\frac{T(x_0)}{\delta_0(x_0)}) )}$.
Moreover, 
 for all $s'\ge s\ge -\log T^*(x)$, the functional
$F_0(w(s),s)$ defined in (\ref{00E3}) satisfies
\begin{eqnarray}\label{LE}
F_0(w(s'),s')-F_0(w(s),s) &=&
- \int_{s}^{s'} e^{2\alpha \tau}\!\int_{\partial B}\!\Big(\partial_s w+\alpha w\Big)^2{\mathrm{d}}\sigma {\mathrm{d}}\tau\nonumber\\
&&+ \frac{\alpha (p-1)}{p+1}\int_{s}^{s'} e^{2\alpha
\tau}\int_{B}|w|^{p+1}{\mathrm{d}}y{\mathrm{d}}\tau.
\end{eqnarray}
Moreover, for all $s\ge -\log T^*(x)$, we have  $F_0(w(s),s)\ge0$,
and
\begin{equation}\label{C0}
F_0(w(s),s)\rightarrow 0
\textrm{ as }s\rightarrow +\infty.
\end{equation}
}

\bigskip

Throughout this section,   we consider $u $   a solution of (\ref{gen}) with
blow-up graph $\Gamma:\{x\mapsto T(x)\}$ and  $x_0$ is a
non-characteristic point.
For any $T_0\in (0, T(x_0)]$ and $x\in \er^N$  such that $|x-x_0|\le
\frac{T_0}{\delta_0(x_0)}$, we will write $w$ instead of
$w_{x,T^*(x)}$ defined in (\ref{scaling}) with $T^*(x)$ given in
(\ref{18dec1}).

\medskip

Let us first note that equation \eqref{C} satisfied by $w$ can be rewritten  in the following form:
\begin{eqnarray}\label{eqw1}
\partial^2_{s}w &=&\frac{1}{\p}\div(\p \grad w-\p (y\cdot \grad w)
y)+(2\varepsilon-2\alpha)
 y\cdot \nabla w
-\frac{2(p+1)}
{(p-1)^2}w+|w|^{p-1}w\nonumber
\\ 
&&-\frac{p+3}{p-1}\partial_{s}w 
-2y\cdot \nabla
\partial_{s}w, \qquad  \forall y\in B\ {\textrm{ and}} \ s\ge -\log T^*(x),
\end{eqnarray}
for any $\e>0,$ where 
\begin{equation}\label{A50}
\p=(1-|y|^2)^{\varepsilon}.
\end{equation}


As we said earlier,  we first introduce two classical energy estimates obtained by 
  multiplying the   equation \eqref{C} by $w\p$ and $\partial_sw\p$,
  before integration.

\subsection{Classical  energy  estimates}

In this subsection, we   introduce  the following 
natural 
functionals:
\begin{eqnarray}
 E_{\e}(w(s))
 &=&\displaystyle\int_B \Big (
\frac{1}{2}(\partial_sw)^2 
+\frac{1}{2}(|\nabla w|^2-(y\cdot\grad w)^2)
+\frac{p+1}{(p-1)^2}w^2
-\frac{|w|^{p+1}}{p+1}\Big )\p{\mathrm{d}}y,\label{En}\qquad
\\
J_{\e}(w(s))&=& -\displaystyle\int_{B}w\partial_s w\p  {\mathrm{d}}y- (\frac{ N}2+\alpha) \displaystyle\int_{B}w^2\p
 {\mathrm{d}}y,\label{In}
\end{eqnarray}
where  $\e >0$.

\medskip

We first estimate the time derivative of $E_{\e}(w(s))$  in the following lemma:
\begin{lem} \label{LE} For all $\e>0$, for all  
  $s'\ge s \geq-\log T^*(x)$, we have 
  \begin{eqnarray}\label{5avril34}
 E_{\e}(w(s'))-E_{\e}(w(s))&=&-2\varepsilon\int_s^{s'}
 \int_B (\partial_{s}w)^2\frac{|y|^2\p}{1-|y|^2}\y{\mathrm{d}}\tau
-2\alpha \int_s^{s'}\int_B (\partial_s w)^2 \p{\mathrm{d}}y{\mathrm{d}}\tau\nonumber\\
&&+(2\varepsilon-2\alpha)\int_s^{s'}\int_B \partial_s w y\cdot \nabla w \p{\mathrm{d}}y{\mathrm{d}}\tau.
\end{eqnarray}
\end{lem}
\noindent {\it Proof}: 
Using \eqref{eqw1} and  integration by parts, we get  for all   $s \geq-\log T^*(x)$,
\begin{equation*}
 \frac{d}{ds}E_{\e}(w(s))=(2\e-2\alpha)\int_B \partial_s w y\cdot \nabla w \p{\mathrm{d}}y
-\frac{p+3}{p-1} \int_B (\partial_s w)^2 \p{\mathrm{d}}y
+\int_B (\partial_s w)^2 \dv (y\p){\mathrm{d}}y.
\end{equation*}
Noticing that 
\begin{equation}\label{div01}
\div (\p y)=N\p-2\e\frac{|y|^2}{1-|y|^2}\p.
\end{equation}
Therefore, exploting  \eqref{div01} and  the definition of $\alpha$ given by \eqref{alpha}, we  infer
\begin{eqnarray}\label{5avril33}
 \frac{d}{ds}E_{\e}(w(s))&=&-2\varepsilon
 \int_B (\partial_{s}w)^2\frac{|y|^2\p}{1-|y|^2}\y
-2\alpha \int_B (\partial_s w)^2 \p{\mathrm{d}}y\nonumber\\
&&+(2\varepsilon-2\alpha)\int_B \partial_s w y\cdot \nabla w \p{\mathrm{d}}y.
\end{eqnarray}
Integrating  \eqref{5avril33} in time yields \eqref{5avril34}.
This ends the proof of Lemma  \ref{LE}.

\Box

Thanks to  Lemma \ref{LE}, we are in position to prove the following:
\begin{lem}\label{ws01}
For all $\e >0$,  for all  $s\ge -\log T^*(x)+1,$ we have
\begin{equation}\label{ws1}
\int_{s}^{s+1}\int_{B}(\partial_s w)^2\frac{\p}{1-|y|^2}\y 
{\mathrm{d}}\tau \le C\int_{s-1}^{s+2}\!\int_{B}\Big(|\nabla  w|^2+(\partial_sw)^2+w^2+|w|^{p+1} \Big) \p {\mathrm{d}}y{\mathrm{d}}\tau. 
\end{equation}
\end{lem}

{\it Proof:}
Let  $s\ge -\log T^*(x)+1$, $s_1=s_1(s)\in [s-1,s]$  
and $s_2=s_2(s)\in [s+1,s+2]$ to be chosen later. From Lemma \ref{LE}, 
by the Cauchy-Schwarz  inequality, and
  the fact that $ab\le a^2+b^2$,  we can write 
\begin{eqnarray}\label{ws2}
\int_{s}^{s+1}\int_{B}(\partial_sw)^2\frac{|y|^2\p}{{1-|y|^2}}\y 
{\mathrm{d}}\tau &\le &\hspace{-0,3cm} C\int_{s_1}^{s_2}\!\int_{B}\Big(|\nabla  w|^2+(\partial_sw)^2+w^2+|w|^{p+1} \Big) \p {\mathrm{d}}y{\mathrm{d}}\tau   \nonumber\\
 &&+|E_{\e }(w(s_2))|+|E_{\e }(w(s_1))|.
\end{eqnarray}
Now, we control all the terms on the right-hand side
of the relation (\ref{ws2}). Note that, by the expression  of $E_{\e }(w(s))$ and using the Cauchy-Schwarz inequality,
we can write, 
\begin{equation}\label{ws3}
|E_{\e }(w(s_i))|\le C \displaystyle\int_{B}\Big(
(w(s_i))^2+|w(s_i)|^{p+1}+
(\partial_s
w(s_i))^2+|\nabla  w(s_i)|^2\Big) \p {\mathrm{d}}y,\qquad i \in \{1,2\}.
\end{equation}
By using the mean value theorem, let us choose  $s_1=s_1(s)\in
[s-1,s]$ such that
\begin{equation}\label{ws4}
\displaystyle\int_{s-1}^{s} \displaystyle\int_{B}\Big(w^2+|w|^{p+1}+(\partial_s
w)^2+|\nabla  w|^2\Big) \p {\mathrm{d}}y{\mathrm{d}}\tau= \int_{B}\Big((w(s_1))^2+|w(s_1)|^{p+1}
+(\partial_s
w(s_1))^2+|\nabla  w(s_1)|^2\Big) \p {\mathrm{d}}y.
\end{equation}
In view of \eqref{ws3},  and (\ref{ws4}) we write
\begin{equation}\label{ws5}
|E_{\e }(w(s_1))|\le C\displaystyle\int_{s-1}^{s} \displaystyle\int_{B}\Big(w^2+|w|^{p+1}+(\partial_s
w)^2+|\nabla  w|^2\Big) \p{\mathrm{d}}y{\mathrm{d}}\tau.
\end{equation}

%
Also, by using the mean value theorem, we choose
$s_2=s_2(s)\in [s+1,s+2]$ such that
\begin{equation}\label{ws7}
\displaystyle\int_{s+1}^{s+2} \displaystyle\int_{B}\Big(w^2+|w|^{p+1}+(\partial_s
w)^2+|\nabla  w|^2\Big) \p {\mathrm{d}}y{\mathrm{d}}\tau= \displaystyle\int_{B}\Big(w(s_2)^2+|w(s_2)|^{p+1}+(\partial_s
w(s_2))^2+|\nabla  w(s_2)|^2\Big) \p {\mathrm{d}}y.
\end{equation}
Hence, by (\ref{ws3}) and \eqref{ws7}, we infer
\begin{equation}\label{ws8}
|E_{\e }(w(s_2))|\le C\displaystyle\int_{s+1}^{s+2} \displaystyle\int_{B}\Big(w^2+|w|^{p+1}+(\partial_s
w)^2+|\nabla  w|^2\Big)\p {\mathrm{d}}y{\mathrm{d}}\tau.
\end{equation}
 By combining 
 (\ref{ws2}), (\ref{ws5}), and (\ref{ws8})  we deduce
\begin{equation}\label{p4bis24}
\int_{s}^{s+1}\int_{B}(\partial_sw)^2\frac{|y|^2\p}{{1-|y|^2}}\y 
{\mathrm{d}}\tau \le
C\displaystyle\int_{s-1}^{s+2} \displaystyle\int_{B}\Big(w^2+|w|^{p+1}+(\partial_s
w)^2+|\nabla  w|^2\Big) \p{\mathrm{d}}y{\mathrm{d}}\tau.
\end{equation}
Finally, by using the fact that $\frac{(\partial_sw)^2}{1-|y|^2}=\frac{|y|^2(\partial_sw)^2}{1-|y|^2}+(\partial_sw)^2$, and \eqref{p4bis24},
 we conclude \eqref{ws1}. This  ends the proof of Lemma \ref{ws01}.
  
\Box

\bigskip

Now,  we evaluate the time derivative of $J_{\e}(w(s))$  in the following lemma:

\begin{lem} \label{JE}  For all $\e>0$, for all  
  $s'>s \geq-\log T^*(x)$, we have 
 \begin{eqnarray}\label{J0}
J_{\e}(w(s'))-J_{\e}(w(s))&=&-2\alpha \int_s^{s'}J_{\e}(w(\tau)){\mathrm{d}}\tau+\int_s^{s'}
G_{\e}(\tau){\mathrm{d}}\tau
+4\e\int_s^{s'}\int_B  \frac{|y|^2}{1-|y|^2}  \partial_sw w\p \y{\mathrm{d}}\tau \nonumber\\
&&+(2\alpha-2\e)\int_s^{s'}\int_B   y\cdot \grad w w\p \y{\mathrm{d}}\tau+2\alpha\int_s^{s'}\displaystyle\int_{B}w\partial_s w\p  {\mathrm{d}}y{\mathrm{d}}\tau\nonumber\\
&&+\alpha (N+2\alpha) \int_s^{s'}\displaystyle\int_{B}w^2\p
 {\mathrm{d}}y{\mathrm{d}}\tau,
\end{eqnarray}
where  
\begin{equation}\label{G1}
G_{\e}(s)=\int_B |\nabla w|^2\p\y
- \int_B  |w|^{p+1}\p\y-\int_B   
(\partial_{s}w+y\cdot \nabla  w)^2\p \y+
 \frac{2 (p+1)}{(p-1)^2} \displaystyle\int_B
w^2\p{\mathrm{d}}y.
\end{equation}
\end{lem}

{\it Proof}: 
Note that  for all   $s \geq-\log T^*(x)$, we get
\begin{eqnarray}\label{j1}
\frac{d}{ds}J_{\e}(w(s))&=&-\int_B   
(\partial_{s}w)^2\p \y-\int_B   w\partial^2_{s}w\p \y -(N+2\alpha )\int_B   w\partial_{s}w\p \y.
\end{eqnarray}
 According to equation $\eqref{eqw1}$, and by integrating by parts, we obtain
\begin{eqnarray}
\frac{d}{ds}J_{\e}(w(s))&=&
-2 \int_B  y\cdot \grad w \partial_sw\p\y-
2 \int_B  \partial_sw w (y\cdot \grad \p)\y-
2 N\int_B   \partial_sw w\p\y\nonumber\\
&&- \int_B  |w|^{p+1}\p\y
-\int_B   
(\partial_{s}w)^2\p \y+\int_B (|\grad w|^2-(y\cdot \grad w)^2)\p\y\nonumber\\
&&+\frac{2p+2}{(p-1)^2}
  \int_B  w^2\p\y+(2\alpha-2\e)\int_B   y\cdot\grad w w\p \y.
 \end{eqnarray}
 Consequently, by
using  $y\cdot \grad \p=-2\e\frac{|y|^2}{1-|y|^2}\p$,  the definition of the functional $G_{\e}(s)$ in
(\ref{G1}), 
we  get
 \begin{eqnarray}\label{JJ0}
\frac{d}{ds}J_{\e}(w(s))&=&-2\alpha J_{\e}(w(s))+
G_{\e}(s)
+4\e\int_B  \frac{|y|^2}{1-|y|^2}  \partial_sw w\p \y \\
&&+(2\alpha-2\e)\int_B   y\cdot \grad w w\p \y+2\alpha\displaystyle\int_{B}w\partial_s w\p  {\mathrm{d}}y+2\alpha (\frac{ N}2+\alpha) \displaystyle\int_{B}w^2\p
 {\mathrm{d}}y,\nonumber
\end{eqnarray}
Integrating  \eqref{JJ0} in time yields \eqref{J0}.
This ends the proof of Lemma  \ref{JE}.

\Box


\subsection{New  energy  estimates} 
In this subsection, we start by introducing the   three crucial   new
  functionals  defined by the following:
\begin{equation}\label{Po}
N_{\e}(w(s))=\int_{B}\Big(y\cdot \nabla w \partial_sw+(y\cdot \nabla w)^2
\Big)\p {\mathrm{d}}y,
\end{equation}
\begin{equation}\label{Po1bis}
 I_{\e}(w(s))= -\displaystyle\int_{B}w(\partial_s w+2y\cdot \grad w ) \frac{\p}{\sqrt{1-|y|^2}}{\mathrm{d}}y- \frac{ N}2 \displaystyle\int_{B}w^2\frac{\p}{\sqrt{1-|y|^2}}
 {\mathrm{d}}y,
\end{equation}
and
\begin{equation}\label{Po1bis1}
 {\cal{L}}_{\e}(w(s))=N_{\frac12+\e}(w(s))+(\frac12+\e){{I}}_{\e}(w(s)),
\end{equation}
where $\e>0.$

\bigskip
 
First, we shall show that the above three functionals are well-defined. 
From    the fact that $(w,\partial_sw) \in  {\cal H}$, where  ${\cal H}$ is defined in \eqref{HH},  a covering technique,  and  the basic inequality
$\big|N_{\e}(w(s))\big|\le C  \int_{B} \big(|\nabla w|^2 +(\partial_sw)^2\big)
 {\mathrm{d}}y$, we conclude that the functional $N_{\e}(w(s))$ is well-defined.
Next,  recall from   \cite{MZajm03} the following Hardy type inequality
\begin{equation*}
 \int_{B}w^2\frac{|y|^2}{{1-|y|^2}}\p\y\leq C\int_{B}|\nabla w|^2(1-|y|^2)\p \y+C\int_{B} w^2\p \y.
\end{equation*}
Hence, using the basic identity  $\frac{w^2}{1-|y|^2}=\frac{|y|^2w^2}{1-|y|^2}+w^2$, and the above  Hardy  inequality,
we conclude that 
\begin{equation}\label{Hardy}
 \int_{B}w^2\frac{\p}{{1-|y|^2}}\y\leq C\int_{B}|\nabla w|^2(1-|y|^2)\p \y+C\int_{B} w^2\p \y.
\end{equation}
By observing
$\big|I_{\e}(w(s))\big|\le C  \int_{B} \big(|\nabla w|^2 +(\partial_sw)^2+w^2\frac{\p}{1-|y|^2}\big)
 {\mathrm{d}}y$, exploiting \eqref{Hardy},
 we deduce that  $I_{\e}(w(s))$ is well defined. 
 Consequently, 
 the above combined with
 the definition of the functional ${\cal{L}}_{\e}(w(s))$ given by \eqref{Po1bis1}, we  easily  deduce that 
${\cal{L}}_{\e}(w(s))$ is also well defined.

\bigskip
%

  Let us denote that the first functional  $N_{\e}(w(s))$ arises from a Pohozaev identity obtained through 
the multiplication of equation \eqref{C} by $(y\cdot \grad w)\p$.  The second functional is obtained by  multipliying the equation \eqref{C} by $w(1-|y|^2)^{-\frac12}\p$.  Finally, we introduce the functional ${\cal{L}}_{\e}(w(s))=N_{\frac12+\e}(w(s))+aI_{\e}(w(s))$, where $a$ is chosen such that the very singular term 
$\int_{B}|\nabla_{\theta}w|^2
\frac{|y|^2\p}{1-|y|^2}\y$
is canceled once we estimate the derivative of the functional    ${\cal{L}}_{\e}(w)$.   This estimate allows us to get control of 
 the time average of the $L^{p+1}$ norm of $w$ with a singular weight.
This is one novelty of our paper. 

\medskip

We begin by estimating the functional ${{N}}_{\e }(w(s))$ in the following lemma:

\begin{lem} \label{L01} 
For all $\e>0$, for all  
  $s'>s \geq-\log T^*(x)$, we have 
 \begin{eqnarray}\label{N0t}
N_{\e}(w(s'))-N_{\e}(w(s))&=&-2\alpha \int_s^{s'}N_{\e}(w(\tau)){\mathrm{d}}\tau
-\varepsilon
\int_s^{s'} \int_{B}|\nabla_{\theta}w|^2
\frac{|y|^2\p}{1-|y|^2}\y{\mathrm{d}}\tau\nonumber\\
&&+\frac{2\e}{p+1}\int_s^{s'}
\int_{B} |w|^{p+1}\frac{\p}{1-|y|^2}\y{\mathrm{d}}\tau\\
&&-\frac{N}2 \int_s^{s'}\int_{B}  (\partial_sw)^2 \p {\mathrm{d}}y{\mathrm{d}}\tau
+\e  \int_s^{s'} \int_{B}  (\partial_sw)^2 \frac{|y|^2\p}{1-|y|^2}
{\mathrm{d}}y{\mathrm{d}}\tau \nonumber\\
&&+\frac{N}2 \int_s^{s'}\int_{B}\Big(|\grad w|^2-(y\cdot \grad  w)^2\Big) \p
{\mathrm{d}}y{\mathrm{d}}\tau - \int_s^{s'}\int_{B}|\grad w|^2 \p
{\mathrm{d}}y{\mathrm{d}}\tau\nonumber\\
&&-\frac{2(p+1)}{(p-1)^2} \int_s^{s'}\int_{B}y\cdot \nabla 
ww\p{\mathrm{d}}y{\mathrm{d}}\tau -\frac{N+2\e}{p+1} \int_s^{s'}\int_{B} |w|^{p+1}\p
 {\mathrm{d}}y{\mathrm{d}}\tau\nonumber\\
 &&-N
 \int_s^{s'}\int_{B}y\cdot \nabla   w \partial_{s}w\p {\mathrm{d}}y{\mathrm{d}}\tau+\e   \int_s^{s'}\int_{B}(y \cdot \nabla   w)^2\p{\mathrm{d}}y{\mathrm{d}}\tau.\nonumber
\end{eqnarray}
\end{lem}
%

%
\bigskip

{\it Proof:} Note that  
for all
 $s\ge -\log T^*(x)$, we have
\begin{equation}\label{N1}
\frac{d}{ds} N_{\e }(w(s))= \int_{B}y\cdot \nabla w (\partial^2_sw+2y\cdot \nabla\partial_s w)
\p {\mathrm{d}}y
+\int_{B}y\cdot \nabla \partial_sw\partial_sw
\p {\mathrm{d}}y.
\end{equation}
By using \eqref{div01}, and  integrating by parts, we infer 
$$\int_{B}y\cdot \nabla \partial_sw\partial_sw
\p {\mathrm{d}}y=-\frac{N}2\int_{B}  (\partial_sw)^2 \p {\mathrm{d}}y
+\e  \int_{B}  (\partial_sw)^2 \frac{|y|^{2}\p}{1-|y|^2}
{\mathrm{d}}y.$$
Combining the above  equality with \eqref{N1}, we may arrive at
\begin{eqnarray}\label{3}
\frac{d}{ds} {{N}}_{\e }(w(s))&=&-\frac{N}2\int_{B}  (\partial_sw)^2 \p {\mathrm{d}}y
+\e  \int_{B}  (\partial_sw)^2 \frac{|y|^{2}\p}{1-|y|^2}
{\mathrm{d}}y \nonumber\\
&& \int_{B}y\cdot \nabla w (\partial^2_sw+2y\cdot \nabla\partial_s w)
\p {\mathrm{d}}y.
\end{eqnarray}
By using \eqref{eqw1},  we have
\begin{eqnarray}\label{N11}
\frac{d}{ds} N_{\e }(w(s))&=&-\frac{N}2\int_{B}  (\partial_sw)^2 \p {\mathrm{d}}y
+\e  \int_{B}  (\partial_sw)^2 \frac{|y|^{2}\p}{1-|y|^2}
{\mathrm{d}}y \nonumber\\
&&+  \int_{B}y\cdot \nabla   w \div(\p \grad w-\p (y\cdot \grad w)
y){\mathrm{d}}y
+(2\e-2\alpha)  \int_{B}(y\cdot \nabla   w)^2\p{\mathrm{d}}y\nonumber\\
&&-\frac{2(p+1)}{(p-1)^2}\int_{B}y\cdot \nabla 
w w\p{\mathrm{d}}y
 +\int_{B}y\cdot \nabla  w |w|^{p-1}w
\p {\mathrm{d}}y\nonumber\\
 &&-\frac{p+3}{p-1}
\int_{B}y\cdot \nabla   w \partial_{s}w\p {\mathrm{d}}y.
\end{eqnarray}
Substituting  \eqref{mport01} into  (\ref{N11}),  and using the fact $\frac{p+3}{p-1}=N+2\alpha$, we have
\begin{eqnarray}\label{N11bis}
\frac{d}{ds} N_{\e }(w(s))&=&-\varepsilon \int_{B}|\nabla_{\theta}w|^2
\frac{|y|^2\p}{1-|y|^2}\y-\frac{N}2\int_{B}  (\partial_sw)^2 \p {\mathrm{d}}y
+\e  \int_{B}  (\partial_sw)^2 \frac{|y|^{2}\p}{1-|y|^2}
{\mathrm{d}}y \nonumber\\
&& +\frac{N}2\int_{B}\Big(|\grad w|^2-(y\cdot\grad  w)^2\Big) \p
{\mathrm{d}}y 
+(\e-2\alpha)  \int_{B}(y\cdot \nabla   w)^2\p{\mathrm{d}}y\nonumber\\
&&-\frac{2(p+1)}{(p-1)^2}\int_{B}y\cdot \nabla 
w w\p{\mathrm{d}}y
 +\int_{B}y\cdot \nabla  w |w|^{p-1}w
\p {\mathrm{d}}y\nonumber\\
 &&-(N+2\alpha )\int_{B}y\cdot \nabla   w \partial_{s}w\p {\mathrm{d}}y
-\int_{B}|\grad w|^2 \p
{\mathrm{d}}y.
\end{eqnarray}
  By  using the expression of $N_{\e }(w(s))$ given by \eqref{Po}, \eqref{div01}, and  integrating by parts 
we conclude 
\begin{eqnarray}\label{N0}
\frac{d}{ds} {{N}}_{\e }(w(s))
&=&-2\alpha  {{N}}_{\e }(w(s))-\varepsilon \int_{B}|\nabla_{\theta}w|^2
\frac{|y|^2\p}{1-|y|^2}\y +\frac{2\e}{p+1}\int_{B} |w|^{p+1}\frac{\p}{1-|y|^2}\y\nonumber\\
&&-\frac{N}2\int_{B}  (\partial_sw)^2 \p {\mathrm{d}}y
+\e  \int_{B}  (\partial_sw)^2 \frac{|y|^2\p}{1-|y|^2}
{\mathrm{d}}y \nonumber\\
&&+\frac{N}2\int_{B}\Big(|\grad w|^2-(y\cdot \grad  w)^2\Big) \p
{\mathrm{d}}y -\int_{B}|\grad w|^2 \p
{\mathrm{d}}y\nonumber\\
&&-\frac{2(p+1)}{(p-1)^2}\int_{B}y\cdot \nabla 
ww\p{\mathrm{d}}y -\frac{N+2\e}{p+1}\int_{B} |w|^{p+1}\p
 {\mathrm{d}}y\nonumber\\
 &&-N
\int_{B}y\cdot \nabla   w \partial_{s}w\p {\mathrm{d}}y+\e  \int_{B}(y \cdot \nabla   w)^2\p{\mathrm{d}}y.
\end{eqnarray}
Integrating  \eqref{N0} in time yields 
 \eqref{N0t}.
This ends the proof of Lemma \ref{L01}. 
\Box

\bigskip

Furthermore,  by estimating  the time derivative of ${I}_{\e }(w(s))$, we conclude  the following lemma:
\begin{lem}\label{L01bis}
For all $\e>0$, for all  
  $s'>s \geq-\log T^*(x)$, we have 
 \begin{eqnarray}\label{V01t}
I_{\e}(w(s'))-I_{\e}(w(s))&=&-2\alpha \int_s^{s'}I_{\e}(w(\tau)){\mathrm{d}}\tau
-\int_s^{s'} \int_{B}|w|^{p+1} \frac{\p}{\sqrt{1-|y|^2}}{\mathrm{d}}y
{\mathrm{d}}\tau
\nonumber\\
&&+\int_s^{s'}\int_{B}|\grad w_{\theta}|^2\pp{\mathrm{d}}y{\mathrm{d}}\tau
-\int_s^{s'}\int_{B}(\partial_sw)^2\pp {\mathrm{d}}y{\mathrm{d}}\tau\nonumber\\
&&+(1-2\e-2\alpha)  \int_s^{s'}\displaystyle\int_{B}
wy\cdot \grad w \frac{\p}{\sqrt{1-|y|^2}}{\mathrm{d}}y{\mathrm{d}}\tau
 \nonumber\\
&&+(\frac{2(p+1)}
{(p-1)^2}-\alpha N)\int_s^{s'}\displaystyle\int_{B} w^2 \frac{\p}{\sqrt{1-|y|^2}}{\mathrm{d}}y{\mathrm{d}}\tau\nonumber\\
&&-2 \int_s^{s'}\displaystyle\int_{B}\partial_s wy\cdot \grad w \frac{\p}{\sqrt{1-|y|^2}}{\mathrm{d}}y{\mathrm{d}}\tau\nonumber\\
&&+\int_s^{s'}\int_{B}|\grad w_r|^2{\rho_{\e+\frac12}}{\mathrm{d}}y{\mathrm{d}}\tau.
\end{eqnarray}
\end{lem}

\bigskip

{\it Proof:} Note that   for all
 $s\ge -\log T^*(x)$, we have
\begin{eqnarray*}
\frac{d}{ds}{I}_{\e }(w(s))&=&- \int_{B}(\partial_sw)^2\pp{\mathrm{d}}y-
\int_{B}w(\partial^2_{s}w
+2y\cdot \grad \partial_sw) 
\pp{\mathrm{d}}y \\
&&-2 \displaystyle\int_{B}\partial_s wy\cdot \grad w \frac{\p}{\sqrt{1-|y|^2}}{\mathrm{d}}y- N \int_{B}w\partial_s w\pp
{\mathrm{d}}y.
\end{eqnarray*}
By using  equation (\ref{eqw1}), 
 we have
\begin{eqnarray}\label{R0}
\frac{d}{ds}{I}_{\e }(w(s))&=&- \int_{B}(\partial_sw)^2\pp {\mathrm{d}}y- \int_{B}\div(\p \grad w-\p (y\cdot \grad w)
y)w\frac1{\sqrt{1-|y|^2}}{\mathrm{d}}y \nonumber\\
&&+2(\alpha-\e) \displaystyle\int_{B}wy\cdot \grad w \frac{\p}{\sqrt{1-|y|^2}}{\mathrm{d}}y
+\frac{2(p+1)}
{(p-1)^2}\displaystyle\int_{B} w^2 \frac{\p}{\sqrt{1-|y|^2}}{\mathrm{d}}y\nonumber\\
&&-\displaystyle\int_{B} |w|^{p+1} \frac{\p}{\sqrt{1-|y|^2}}{\mathrm{d}}y+(\frac{p+3}{p-1}-N) \displaystyle\int_{B} \partial_sww \frac{\p}{\sqrt{1-|y|^2}}{\mathrm{d}}y\nonumber\\
&&-2 \displaystyle\int_{B}\partial_s wy\cdot \grad w \frac{\p}{\sqrt{1-|y|^2}}{\mathrm{d}}y.
\end{eqnarray}
Substituting  \eqref{mport1} into  (\ref{R0}),  and using the fact $\frac{p+3}{p-1}-N=2\alpha$, we have
\begin{eqnarray}\label{V01}
\frac{d}{ds}{I}_{\e }(w(s))&=&-2\alpha {I}_{\e }(w(s)) - \int_{B}(\partial_sw)^2\pp {\mathrm{d}}y
+\int_{B}|\grad w_{\theta}|^2\pp{\mathrm{d}}y
 \nonumber\\
&&-\displaystyle\int_{B} |w|^{p+1} \frac{\p}{\sqrt{1-|y|^2}}{\mathrm{d}}y+(1-2\e-2\alpha) \displaystyle\int_{B}
wy\cdot \grad w \frac{\p}{\sqrt{1-|y|^2}}{\mathrm{d}}y\nonumber\\
&&+(\frac{2(p+1)}
{(p-1)^2}-\alpha N)\displaystyle\int_{B} w^2 \frac{\p}{\sqrt{1-|y|^2}}{\mathrm{d}}y
-2 \displaystyle\int_{B}\partial_s wy\cdot \grad w \frac{\p}{\sqrt{1-|y|^2}}{\mathrm{d}}y\nonumber\\
&&+\int_{B}|\grad w_r|^2{\rho_{\e+\frac12}}{\mathrm{d}}y.
\end{eqnarray}
Integrating  \eqref{V01} in time yields \eqref{V01t}.
This ends the proof of Lemma \ref{L01bis}. 
\Box

\medskip

Combining the above inequalities with Lemmas   \ref{L01} and \ref{L01bis}, we consequently reach the
 following lemma:
\begin{lem}\label{L01bisbis}
For all $\e>0$, for all  $s'> s\ge -\log T^*(x)$, we have 
\begin{eqnarray}\label{R02bis}
 {\cal{L}}_{\e }(w(s'))
&=& {\cal{L}}_{\e }(w(s))-2\alpha \int_s^{s'} 
{\cal{L}}_{\e }(w(\tau)){\mathrm{d}}\tau
 -\frac{(1+2\e)(p-1)}{2(p+1)} \int_s^{s'} \int_{B} |w|^{p+1}\frac{\p}{\sqrt{1-|y|^2}}\y{\mathrm{d}}\tau
 \nonumber\\
&&+\frac{N}2 \int_s^{s'} \int_{B}\Big(|\grad w|^2-(y\cdot\grad  w)^2\Big) \rho_{\frac12+\e}
{\mathrm{d}}y{\mathrm{d}}\tau -(\frac12-\e) \int_s^{s'} \int_{B}|\grad w|^2 \rho_{\frac12+\e}
{\mathrm{d}}y{\mathrm{d}}\tau \nonumber\\
&&-\frac{2(p+1)}{(p-1)^2} \int_s^{s'} \int_{B}y\cdot\nabla 
ww\rho_{\frac12+\e}{\mathrm{d}}y{\mathrm{d}}\tau -\frac{N+1+2\e}{p+1} \int_s^{s'} \int_{B} |w|^{p+1}\rho_{\frac12+\e}
 {\mathrm{d}}y{\mathrm{d}}\tau\nonumber\\
 &&-N \int_s^{s'} 
\int_{B}y\cdot\nabla   w \partial_{s}w\rho_{\frac12+\e} {\mathrm{d}}y{\mathrm{d}}\tau+(\frac12+\e)  \int_s^{s'}  \int_{B}(y\cdot \nabla   w)^2\rho_{\frac12+\e}
{\mathrm{d}}y{\mathrm{d}}\tau\nonumber\\
&&- (\frac{N+1}2+\e) \int_s^{s'} 
\int_{B}(\partial_sw)^2\rho_{\frac12+\e} {\mathrm{d}}y{\mathrm{d}}\tau\nonumber\\
&&+(1-2\e-2\alpha ) (\frac12+\e) \int_s^{s'} \displaystyle\int_{B}
wy\cdot \grad w \frac{\p}{\sqrt{1-|y|^2}}{\mathrm{d}}y{\mathrm{d}}\tau\nonumber\\
&&+(\frac{2(p+1)}
{(p-1)^2}-\alpha N)(\frac12+\e) \int_s^{s'} \displaystyle\int_{B} w^2 \frac{\p}{\sqrt{1-|y|^2}}{\mathrm{d}}y{\mathrm{d}}\tau\nonumber\\
&&-(1+2\e) \int_s^{s'}  \displaystyle\int_{B}\partial_s wy\cdot \grad w \frac{\p}{\sqrt{1-|y|^2}}{\mathrm{d}}y{\mathrm{d}}\tau.
\end{eqnarray}
\end{lem}

{\it Proof:} Note that  
for all $\e>0$, for all  $s'> s\ge -\log T^*(x)$, we have 
\begin{eqnarray}\label{N0t1}
N_{\frac12+\e}(w(s'))&=&N_{\frac12+\e}(w(s))-2\alpha \int_s^{s'}N_{\frac12+\e}(w(\tau)){\mathrm{d}}\tau
-(\frac12+\varepsilon)
\int_s^{s'} \int_{B}|\nabla_{\theta}w|^2
\frac{|y|^2\p}{\sqrt{1-|y|^2}}\y{\mathrm{d}}\tau\nonumber\\
&&+\frac{1+2\e}{p+1}\int_s^{s'}
\int_{B} |w|^{p+1}\frac{\p}{\sqrt{1-|y|^2}}\y{\mathrm{d}}\tau\\
&&-\frac{N}2 \int_s^{s'}\int_{B}  (\partial_sw)^2 \rho_{\frac12+\e} {\mathrm{d}}y{\mathrm{d}}\tau
+({\frac12+\e})  \int_s^{s'} \int_{B}  (\partial_sw)^2 \frac{|y|^2\p}{\sqrt{1-|y|^2}}
{\mathrm{d}}y{\mathrm{d}}\tau \nonumber\\
&&+\frac{N}2 \int_s^{s'}\int_{B}\Big(|\grad w|^2-(y\cdot \grad  w)^2\Big) \rho_{\frac12+\e}
{\mathrm{d}}y{\mathrm{d}}\tau - \int_s^{s'}\int_{B}|\grad w|^2 \rho_{\frac12+\e}
{\mathrm{d}}y{\mathrm{d}}\tau\nonumber\\
&&-\frac{2(p+1)}{(p-1)^2} \int_s^{s'}\int_{B}y\cdot \nabla 
ww\rho_{\frac12+\e}{\mathrm{d}}y{\mathrm{d}}\tau -\frac{N+1+2\e}{p+1} \int_s^{s'}\int_{B} |w|^{p+1}\rho_{\frac12+\e}
 {\mathrm{d}}y{\mathrm{d}}\tau\nonumber\\
 &&-N
 \int_s^{s'}\int_{B}y\cdot \nabla   w \partial_{s}w\rho_{\frac12+\e}{\mathrm{d}}y{\mathrm{d}}\tau+(\frac12+\e)   \int_s^{s'}\int_{B}(y \cdot \nabla   w)^2\rho_{\frac12+\e}{\mathrm{d}}y{\mathrm{d}}\tau.\nonumber
\end{eqnarray}

Thanks to
\eqref{wr2}, the identities \eqref{V01t}, and \eqref{N0t1},  we get
 yields \eqref{R02bis}.
This ends the proof of Lemma \ref{L01bisbis}.

\Box

\subsection{Proof of Proposition \ref{P01m}':}

Now, we are in position to state and prove Proposition \ref{P01m}', which is a uniform
 version of  Proposition \ref{P01m}  for $x$ near $x_0$.\\

\noindent {\bf{Proposition}\ref{P01m}'}\\
Let $\e >0$. Let  $u $    a solution of ({\ref{gen}})
with blow-up graph $\Gamma :\{x\mapsto T(x)\}$ and  $x_0$ is a non
characteristic point. Then 
for all $T_0\in (0, T(x_0)]$, for all
 $s\ge 2-\log (T_0)$ and $x\in \er^N$, where $|x-x_0|\le \frac{e^{-s}}{\delta_0(x_0)}$,
we have

\begin{equation}\label{w1}
\int_{s}^{s+1}\int_{B}|w|^{p+1}\frac{\p}{\sqrt{1-|y|^2}}\y 
{\mathrm{d}}\tau \le C\int_{s-2}^{s+3}\!\int_{B}\Big(|\nabla  w|^2+(\partial_sw)^2+w^2+|w|^{p+1} \Big) \p {\mathrm{d}}y{\mathrm{d}}\tau,
\end{equation}
where $w=w_{x,T^*(x)}$,  where $T^*(x)$ is defined \eqref{18dec1} and
where the constant $C$ depends only on $N$, and   $p$.

\bigskip

{\it Proof:}
Let  $s\ge 2-\log T^*(x)$, $s_3=s_3(s)\in [s-1,s]$  
and $s_4=s_4(s)\in [s+1,s+2]$ to be chosen later. From Lemma \ref{L01bisbis}, 
by the Cauchy-Schwarz  inequality, and
  the fact that $ab\le a^2+b^2$,  we can write 
\begin{eqnarray}\label{p4bis}
\int_{s}^{s+1}\int_{B}\Big||w|^{p+1}\frac{\p}{\sqrt{1-|y|^2}}\y 
{\mathrm{d}}\tau &\le &\hspace{-0,2cm} C\int_{s_3}^{s_4}\!\int_{B}\Big(|\nabla  w|^2+\frac{(\partial_sw)^2}{1-|y|^2}+\frac{w^2}{1-|y|^2}+|w|^{p+1} \Big) \p {\mathrm{d}}y{\mathrm{d}}\tau   \nonumber\\
 &&\hspace{-0,35cm}+C \Big|\int_{s_3}^{s_4} 
{\cal{L}}_{\e }(w(\tau)){\mathrm{d}}\tau\Big| 
 +C|{\cal{L}}_{\e }(w(s_4))|+C|
{\cal{L}}_{\e }(w(s_3))|.\nonumber\\
\end{eqnarray}
Now, we control all the terms on the right-hand side
of the relation (\ref{p4bis}). Note that, by the expression  of ${\cal L}_{\e }(w(s))$ and using the Cauchy-Schwarz inequality,
we can write
\begin{equation}\label{I01}
|{\cal{L}}_{\e }(w(s))|\le C \displaystyle\int_{B}\Big(
\frac{(w(s))^2}{1-|y|^2}+
(\partial_s
w(s))^2+|\nabla  w(s)|^2\Big) \p {\mathrm{d}}y,\qquad \forall s\in [s_3,s_4].
\end{equation}
By using the mean value theorem, let us choose  $s_3=s_3(s)\in
[s-1,s]$ such that
\begin{equation}\label{q04}
\displaystyle\int_{s-1}^{s} \displaystyle\int_{B}\Big(\frac{w^2}{1-|y|^2}+(\partial_s
w)^2+|\nabla  w|^2\Big)\p {\mathrm{d}}y{\mathrm{d}}\tau= \displaystyle\int_{B}\Big(\frac{(w(s_3))^2}{1-|y|^2}+(\partial_s
w(s_3))^2+|\nabla  w(s_3)|^2\Big)\p {\mathrm{d}}y.
\end{equation}
For the estimate of $
|{\cal{L}}_{\e }(w(s_3))|$, it follows from
  (\ref{I01}) and (\ref{q04}) that
\begin{equation}\label{I11110}
|{\cal{L}}_{\e }(w(s_3))|\le C\displaystyle\int_{s-1}^{s} \displaystyle\int_{B}\Big(\frac{w^2}{1-|y|^2}+(\partial_s
w)^2+|\nabla  w|^2\Big)\p {\mathrm{d}}y{\mathrm{d}}\tau.
\end{equation}
Similarly, 
by using the mean value theorem, we choose
$s_4=s_4(s)\in [s+1,s+2]$ such that
\begin{equation}\label{q440}
\displaystyle\int_{s+1}^{s+2} \displaystyle\int_{B}\Big(\frac{w^2}{1-|y|^2}+(\partial_s
w)^2+|\nabla  w|^2\Big)\p {\mathrm{d}}y{\mathrm{d}}\tau= \displaystyle\int_{B}\Big(\frac{(w(s_4))^2}{1-|y|^2}+(\partial_s
w(s_4))^2+|\nabla  w(s_4)|^2\Big) \p{\mathrm{d}}y.
\end{equation}
Hence, by (\ref{I01}) and \eqref{q440}, we infer
\begin{equation}\label{I2n}
|{\cal{L}}_{\e }(w(s_4))|\le C\displaystyle\int_{s+1}^{s+2} \displaystyle\int_{B}\Big(\frac{w^2}{1-|y|^2}+(\partial_s
w)^2+|\nabla  w|^2\Big) \p {\mathrm{d}}y{\mathrm{d}}\tau.
\end{equation}
From \eqref{I01}, we get
\begin{equation}\label{p4bis0}
\int_{s_3}^{s_4}
{\cal{L}}_{\e }(w(\tau)){\mathrm{d}}\tau\le 
\ C\int_{s_3}^{s_4}\!\int_{B}\Big(|\nabla  w|^2+\frac{(\partial_sw)^2}{1-|y|^2}+\frac{w^2}{1-|y|^2}+|w|^{p+1} \Big) \p {\mathrm{d}}y{\mathrm{d}}\tau.  
\end{equation}
Gathering  (\ref{p4bis}), (\ref{I11110}), \eqref{I2n}, and (\ref{p4bis0})  we obtain
\begin{equation}\label{p4}
\int_{s}^{s+1}\int_{B}\Big||w|^{p+1}\frac{\p}{\sqrt{1-|y|^2}}\y 
{\mathrm{d}}\tau \le C\int_{s-1}^{s+2}\!\int_{B}\Big(|\nabla  w|^2+\frac{(\partial_sw)^2}{1-|y|^2}+\frac{w^2}{1-|y|^2}+|w|^{p+1} \Big) \p {\mathrm{d}}y{\mathrm{d}}\tau. 
\end{equation}
Now,  then thanks to   Lemma \ref{ws01},  we have  for all  $s\ge -\log T^*(x)+2$ 
\begin{equation}\label{ws11}
\int_{s-1}^{s+2}\int_{B}(\partial_s w)^2\frac{\p}{1-|y|^2}\y 
{\mathrm{d}}\tau \le C\int_{s-2}^{s+3}\!\int_{B}\Big(|\nabla  w|^2+(\partial_sw)^2+w^2+|w|^{p+1} \Big) \p {\mathrm{d}}y{\mathrm{d}}\tau. 
\end{equation}
By exploiting  \eqref{p4}, \eqref{ws11}, and \eqref{Hardy}, we easily deduce
the desired estimate $\eqref{w1}$. This ends the proof of Proposition \ref{P01m}'.

\Box

\section{Improvments about  the $\dot{H}^{1}(\R^N\times [-\log(s_0),\infty))$ norm of $w$ 
  }\label{sec3}


 Consider $u $   a solution of (\ref{gen}) with
blow-up graph $\Gamma:\{x\mapsto T(x)\}$ and  $x_0$ is a non
characteristic point.
Let   $T_0\in (0, T(x_0)]$, for all 
$x\in \er^N$  such that $|x-x_0|\le \frac{T_0}{\delta_0(x_0)}$, then we write $w$ instead of $w_{x,T^*(x)}$ defined in (\ref{scaling}) with $T^*(x)$ given in  (\ref{18dec1}).

\par

First, although the computations in  this section can be realized for all  $\e \in (\frac12,1)$ 
for simplicity, we shall assume 
from now to set  $\e=\e_0=\frac35$  and ignore the dependence on $\e$ for all the functions 
 that will be subsequently used in this section. For
example, the function $\rho_{\e_0}$
 will be simply denoted by $\rho$ for simplicity, and the same for all the other
functions (including the constants).

This section is divided into two parts:
\begin{itemize}
\item  The  subsection \ref{241}, is devoted to deriving an 
improvement about the $\dot{H}^{1}(\R^N\times [-\log s_0,\infty))$ norm of $w$, by using the crutial  result showed in Proposition\ref{P01m}'.
\item In subsection \ref{32},
we shall generalize the result obtained in  subsection \ref{241}   under some additional hypotheses.  
Clearly, this finding is useful in the induction argument.
\end{itemize}

\par

\subsection{Proof of Proposition \ref{improv}}\label{241}
As explained before,  
in this subsection, our aim is to prove   the following:
\begin{pro}\label{improv}
 For  all $s'\ge s\ge \max(-\log T^*(x),1)$, we have 
\begin{equation}\label{imp1}
 \int_{s}^{s'}\tau^{-\frac{17}{18}}e^{2\alpha \tau}
   \int_{B}\Big( (\partial_s w)^2+ |\grad w|^2\Big) {\mathrm{d}}y
 {\mathrm{d}}\tau   \le K _4.
 \end{equation}  
\end{pro}

First, we 
combine  \eqref{F0}, \eqref{FF1}, and \eqref{w01}, in the particular case $\e=\e_0=\frac35$, to  deduce the  following crucial estimate:
\begin{cor}\label{P02}
For all $s>\tilde s_0 $,  we have the following inequality
\begin{equation}\label{w2}
\int^{s}_{\tilde s_0}s^{-\frac{17}{18}} e^{2\alpha s}
\int_{B}|w|^{p+1}\frac{\rho}{1-|y|^2} {\mathrm{d}}y{\mathrm{d}}s \leq  K_5,
\end{equation}
where $\tilde s_0=2+\max (-\log T^*(x),1)$.
\end{cor}

\noindent{\it Proof 
}
 By interpolation, we write  for all  $s\ge \tilde s_0$
\begin{multline}\label{N2v1}
\int_{s}^{s+1}e^{2\alpha \tau}\int_{B}| w(\tau )|^{p+1}
\frac{\rho}{1-|y|^2} {\mathrm{d}}y{\mathrm{d}}\tau
\le C\Big(
\int_{s}^{s+1}e^{2\alpha \tau} \int_{B}| w(\tau )|^{p+1}
\frac{1}
{(1-|y|^2)^{\frac9{20}}}
{\mathrm{d}}y{\mathrm{d}}\tau
\Big)^{\frac89} \\
\Big( \int_{s}^{s+1}e^{2\alpha \tau} \int_{B}| w(\tau )|^{p+1} {\mathrm{d}}y{\mathrm{d}}\tau
\Big)^{\frac19}.
\end{multline}
Combining 
 \eqref{w01}, \eqref{cor3}, \eqref{F0}, and \eqref{FF1}, we infer for all $s\ge  \tilde s_0$
\begin{equation}
\int_{s}^{s+1} e^{2\alpha \tau}\int_{B}| w(\tau )|^{p+1}
\frac{\rho}
{{1-|y|^2}}
 {\mathrm{d}}y{\mathrm{d}}\tau
\le C
\Big( \int_{s}^{s+1}e^{2\alpha \tau}\int_{B}| w(\tau )|^{p+1} {\mathrm{d}}y{\mathrm{d}}\tau
\Big)^{\frac19}.
\end{equation}
By using the fact that   $\tau^{-\frac{17}{18}}\le s^{-\frac{17}{18}} $,  for all $\tau \in [s,s+1],$ we get
\begin{eqnarray}\label{N2v1bis}
\int_{s}^{s+1}\tau^{-\frac{17}{18}} e^{2\alpha \tau}\int_{B}| w(\tau )|^{p+1}
\frac{\rho}{{1-|y|^2}}
 {\mathrm{d}}y{\mathrm{d}}\tau
\le Cs^{-\frac{17}{18}} \qquad \qquad\qquad\qquad\qquad\\
\Big( \int_{s}^{s+1}e^{2\alpha \tau}\int_{B}| w(\tau )|^{p+1} {\mathrm{d}}y{\mathrm{d}}\tau
\Big)^{\frac19}, \qquad  \forall s\ge   \tilde s_0.\nonumber
\end{eqnarray}
Thanks to  the basic inequality  $|X|^{1-\eta}|Y|^{\eta}\le C|X|+C |Y|$,  for all $ X,Y\in \R,\eta\in (0,1)$,  we conclude that for all $s\ge  \tilde s_0$
\begin{multline}\label{xv1}
\int_{s}^{s+1}\tau^{-\frac{17}{18}} e^{2\alpha \tau}\int_{B}| w(\tau )|^{p+1}
\frac{\rho}
{{1-|y|^2}}
 {\mathrm{d}}y{\mathrm{d}}\tau
\le  Cs^{
-\frac{17}{16}}\\
+C
\Big( \int_{s}^{s+1}e^{2\alpha \tau}\int_{B}| w(\tau )|^{p+1} {\mathrm{d}}y{\mathrm{d}}\tau
\Big).
\end{multline}
Therefore,  the estimate  \eqref{xv1} implies 
\begin{equation}\label{xv4}
\int_{\tilde s_0}^{\infty}s^{-\frac{17}{18}} e^{2\alpha s}\int_{B}| w(s )|^{p+1}
\frac{\rho}
{{1-|y|^2}}
 {\mathrm{d}}y{\mathrm{d}}s
\le  C+C
\Big( \int_{\tilde s_0}^{\infty}e^{2\alpha s}\int_{B}| w(s )|^{p+1} {\mathrm{d}}y{\mathrm{d}}s
\Big).
\end{equation}
Combining with \eqref{F0}  and \eqref{xv4} leads to  \eqref{w2}.
 This concludes the proof of  Corollary \ref {P02}.
\Box
%

 \bigskip

Now, let us define, for all $s\ge -\log T^*(x)$, 
the following 
 functional:
\begin{eqnarray}\label{Mdef}
{\cal M}(w(s))
&=&E_{\e_0}(w(s))+N_{\e_0}(w(s))-(\frac2{p-1}+\frac2{5}) J_{\e_0}(w(s))\\
 &&+\frac65(\frac2{p-1}+\frac25 )\int_Bw^2  \frac{|y|^2\rho}{1-|y|^2}   \y 
  +(\frac{2}{p-1}+\frac2{5} )\alpha\displaystyle\int_{B} w^2\rho  {\mathrm{d}}y, \nonumber
  \end{eqnarray}

Now, let us first control  the  functional ${\cal M}(w(s))$ in the following lemma:
\begin{lem}\label{Aa1}
For  all  $s'> s\ge -\log T^*(x)$, we have 
\begin{eqnarray}\label{TA1}
{\cal{M}}(w(s'))
&=& {\cal{M}}(w(s))-2\alpha \int_s^{s'} \!\!
{\cal{M}}(w(\tau)){\mathrm{d}}\tau
-\frac35\int_s^{s'} \!\! \int_{B}|\nabla_{\theta}w|^2
\frac{|y|^2\rho}{1-|y|^2}\y{\mathrm{d}}\tau \nonumber\\
&& -\frac9{10}\int_s^{s'} \!\!\int_{B}|\grad w|^2 \rho
{\mathrm{d}}y{\mathrm{d}}\tau
+\frac{6}{5p+5}\int_s^{s'} \!\!\int_{B} |w|^{p+1}\frac{\rho}{1-|y|^2}\y{\mathrm{d}}\tau\nonumber\\
&& +\frac{17-5p}{10p+10}\int_s^{s'} \!\!\displaystyle\int_B 
|w|^{p+1}\rho{\mathrm{d}}y{\mathrm{d}}\tau+\int_s^{s'} \!\Sigma (\tau){\mathrm{d}}\tau,
\end{eqnarray}
where
\begin{equation}\label{sigmasigma}
\int_s^{s'}\! \Sigma (\tau){\mathrm{d}}\tau  +\frac1{20} \int_s^{s'}\!\! \int_B  (\partial_{s}w)^2\rho \y{\mathrm{d}}\tau \le  \frac{4}{5} \int_s^{s'}\!\! \int_{B}|\grad w|^2\rho \y{\mathrm{d}}\tau 
+C_1\int_s^{s'} \!\!\int_{B} w^2\rho  \y{\mathrm{d}}\tau,
\end{equation}
for some $C_1=C_1(p,N)>0$.
\end{lem} 
{\it Proof:} Note that   by using  Lemma \ref{LE}, and the definition of the  $E_{\e_0}(w(s))$  given by \eqref{En}, we can write 
 \begin{eqnarray}\label{5avril3}
 \frac{d}{ds}E_{\e_0}(w(s))&=&-2\alpha E_{\e_0}(w(s))+\alpha G_{\e_0}(s)
+\frac{\alpha (p-1)}{p+1}\displaystyle\int_B 
|w|^{p+1}\rho{\mathrm{d}}y\nonumber\\
&& -\frac65
 \int_B (\partial_{s}w)^2\frac{|y|^2\rho}{1-|y|^2}\y
+\frac65\int_B \partial_s w y\cdot \nabla w \rho{\mathrm{d}}y,
\end{eqnarray}
where $G_{\e_0}(s)$ is given by \eqref{G1}.  Gathering  \eqref{5avril3}, Lemmas  \ref{JE}, and \ref{L01}, we infer
\begin{eqnarray}\label{M00}
\frac{d}{ds}{\cal M}(w(s))&=&   - 2\alpha {\cal M}(w(s))
-\frac35 \int_{B}|\nabla_{\theta}w|^2
\frac{|y|^2\rho}{1-|y|^2}\y  -\frac9{10}\int_{B}|\grad w|^2 \rho
{\mathrm{d}}y\\
&&+\frac{6}{5p+5}\int_{B} |w|^{p+1}\frac{\rho}{1-|y|^2}\y
 +\frac{17-5p}{10p+10}\displaystyle\int_B 
|w|^{p+1}\rho{\mathrm{d}}y+\Sigma (s),\nonumber
\end{eqnarray}
where  $\Sigma(s):=\Sigma_{1}(s)+\Sigma_{2}(s)+\Sigma_{3}(s)+\Sigma_{4}(s),$ and 
\begin{eqnarray*}
\Sigma_{1}(s): &=&-\frac35  \int_{B}  (\partial_sw)^2 \frac{|y|^2\rho}{1-|y|^2}
{\mathrm{d}}y+\frac35  \int_{B}(y \cdot \nabla   w)^2\rho{\mathrm{d}}y
+\frac65\int_B \partial_s w y\cdot \nabla w \rho {\mathrm{d}}y,\nonumber\\
\Sigma_{2}(s) :&=&-\frac1{10}\int_B   (\partial_{s}w)^2\rho \y-\frac1{10}\int_B   (y\cdot \nabla  w)^2\rho \y-\frac1{5}\int_B   
\partial_{s}wy\cdot \nabla  w\rho \y,\nonumber\\
\Sigma_{3}(s) :&=&-\Big(\frac{2p+2}{(p-1)^2}+(\frac2{p-1}+\frac25)(2\alpha-\frac65)\Big)\int_B   y\cdot \grad w w\rho \y \nonumber\\
&&+\Big(\frac{ (p+1)}{5(p-1)^2}-2\alpha (\frac2{p-1}+\frac25) (\frac2{p-1}+\frac{ 1}2+\alpha) \Big)\displaystyle\int_{B}w^2\rho {\mathrm{d}}y\\
\Sigma_{4}(s) :&=&
 -\frac{12\alpha}5(\frac2{p-1}+\frac25 )\int_B w^2 \frac{|y|^2\rho}{1-|y|^2}   \y .
\end{eqnarray*}
Clearly by integrating  \eqref{M00} in time yields \eqref{TA1}. It remains to prove  \eqref{sigmasigma}.

 \bigskip
 
 We are going now to estimate the different terms of $\Sigma(s)$.
By using  the fact that  $y\cdot \grad w=y\cdot \grad_r w$ and the basic inequality  $2ab\le a^2+b^2$, to write for all 
$s\ge -\log T^*(x)$
\begin{equation}\label{sigma1}
2\int_B \partial_s w y\cdot \nabla w \rho{\mathrm{d}}y\le \int_{B}(1-|y|^2)|\nabla_r   w|^2\rho {\mathrm{d}}y+
 \int_B (\partial_{s}w)^2\frac{|y|^2\rho}{1-|y|^2}\y.
\end{equation}
Thanks to the fact  $(1-|y|^2)|\nabla_r   w|^2\le |\grad w|^2-(y.\nabla   w)^2$,  and  the inequality \eqref{sigma1}, we
infer
\begin{equation}\label{5dec5}
\Sigma_1(s)\le \frac35  \int_{B}|\nabla   w|^2\rho {\mathrm{d}}y.
\end{equation}
Similarly, by using the basic inequality $2ab\le \frac12a^2+2b^2$, we write
\begin{eqnarray}\label{later1}
\Sigma_{2}(s)&\le &-\frac1{20} \int_B  (\partial_{s}w)^2\rho \y+\frac1{10} \int_B  (y\cdot \nabla  w)^2\rho \y\nonumber\\
&\le &-\frac1{20} \int_B  (\partial_{s}w)^2\rho \y+\frac1{10} \int_B  |\nabla  w|^2\rho \y.
\end{eqnarray}
Let $A_0=
\Big|\frac{2p+2}{(p-1)^2}+(\frac2{p-1}+\frac25)(2\alpha-\frac65)\Big|$.
By using the inequality   $|y\cdot \grad w| |w|\le \frac{1}{20A_0}|\grad w|^2+C w^2,$ 
we obtain
\begin{equation}\label{5dec1}
\Sigma_{3}(s)  \leq  \frac{1}{20} \int_{B}|\grad w|^2\rho \y
+C\int_{B} w^2\rho  \y.
\end{equation}
We would like now to find an estimate from the term $\Sigma_4(s)$.
Firstly, we use
 the identity $\frac{|y|^2\rho}{1-|y|^2}=-\frac56y\cdot \grad \rho$, to write
\begin{equation}\label{Hg}
 \int_B w^2 \frac{|y|^2\rho}{1-|y|^2}   \y =-\frac56  \int_B w^2 y\cdot \grad \rho\y.
\end{equation}
Hence, by using  \eqref{Hg}, and integrating by parts, we get  
\begin{equation*}
 \int_B w^2 \frac{|y|^2\rho}{1-|y|^2}   \y =\frac{5N}6   \int_B w^2  \rho\y +\frac53 \int_B  y\cdot \grad w w \rho\y.
\end{equation*}
Therefore,  a similar way to the treatment of $\Sigma_3(s)$ leads to 
\begin{equation}\label{5dec123}
\Sigma_{4}(s)  \leq  \frac{1}{20} \int_{B}|\grad w|^2\rho \y
+C\int_{B} w^2\rho  \y.
\end{equation}
Gathering   \eqref{5dec5}, \eqref{later1}, \eqref{5dec1}, and \eqref{5dec123},  we conclude
\begin{equation}\label{5dec123bis}
\Sigma(s)  \leq -\frac1{20} \int_B  (\partial_{s}w)^2\rho \y+ \frac{4}{5} \int_{B}|\grad w|^2\rho \y
+C_2\int_{B} w^2\rho  \y,
\end{equation}
for some $C_2=C_2(p,N)>0$.
A simple integration between $s$ and $s'$ ensures the result  \eqref{sigmasigma}. This concludes the proof 
  of Lemma \ref{Aa1}. \Box \\

\par

Now, we shall  introduce the following functionals:
\begin{eqnarray}
{\cal U}_{1}(w(s),s)&=&\int_{s}^{\infty}\!\tau^{-\frac{17}{18}}e^{2\alpha \tau}\int_{B}\! |w|^{p+1}
\frac{ \rho}{1-|y|^2}  {\mathrm{d}}y {\mathrm{d}}\tau+\int_{s}^{\infty}\!\tau^{-\frac{17}{18}}e^{2\alpha \tau}\int_{B}\! w^2
 \rho {\mathrm{d}}y {\mathrm{d}}\tau,\qquad \label{U}\\
{\cal F}_{1}(w(s),s)&= &s^{-\frac{17}{18}}e^{2\alpha s} {\cal M}(w(s)) +\sigma_1 {\cal U}_1(w(s),s),\label{FF}
\end{eqnarray}
for all $s\ge\tilde s_0,
$ where $\sigma_1$ is a positive constant chosen later.
 \begin{rem}
Let us point out   that, the  functional  ${\cal U}_{1}(w(s),s)$ given by \eqref{U} is well defined.
This is a consequence of 
our previous result  given by \eqref{cor3}  and  the Corollary \ref{P02}. Furthermore, by using the estimates  \eqref{cor3} and \eqref{w2} we  can conclude that  
\begin{eqnarray}\label{LL}
{\cal U}_{1 }(w(s),s) \rightarrow 0 ,\ \ \textrm{ as}\ \ \ s\rightarrow +\infty.
\end{eqnarray}
\end{rem}

By using 
Lemma \ref{Aa1}, we will  choose $\sigma_1>0$, such that  the function ${\cal F}_{1}(w(s),s)$
is  a decreasing functional for the equation
$\eqref{eqw1}$. More precisely, we now  state and prove the following proposition:
\begin{pro}\label{121} There exists $S_1(N,p)>1$, such that 
for  all $s'\ge s\ge  \tilde s_1:=\max(2-\log T^*(x),S_1)$, we have  the following 
\begin{equation}\label{111}
{\cal F}_{1}(w(s'),s') +
\frac1{40} \int_{s}^{s'}\tau^{-\frac{17}{18}}e^{2\alpha \tau}
   \int_{B}\Big( (\partial_s w)^2+ |\grad w|^2\Big) \rho {\mathrm{d}}y
 {\mathrm{d}}\tau   \le {\cal F}_{1}(w(s),s) .
 \end{equation}  
\end{pro}

%

{\it Proof }: A simple calculation yields,  for  all $ s\ge \max(-\log T^*(x),1)$, we have 
\begin{eqnarray}\label{TA10}
 \frac{d}{ds}\Big(  s^{-\frac{17}{18}}e^{2\alpha s}{\cal M}(w(s))\Big)&=&
  s^{-\frac{17}{18}}e^{2\alpha s} 
\Big( \frac{d}{ds}{\cal M}(w(s))  
   + 2\alpha  {\cal M}(w(s))\Big)\nonumber \\
   &&-\frac{17}{18}  s^{-\frac{35}{18}}e^{2\alpha s} {\cal M}(w(s)).
\end{eqnarray}
Therefore, 
by exploiting the  definition of ${\cal F}_{1}(w(s),s)$ given in \eqref{U}, 
estimates \eqref{TA1} and \eqref{TA10}  lead to  
\begin{eqnarray}\label{2411a}
\frac{d}{ds}{\cal F}_{1}(w(s),s)&=&   
-\frac35   s^{-\frac{17}{18}}e^{2\alpha s}\int_{B}|\nabla_{\theta}w|^2
\frac{|y|^2\rho}{1-|y|^2}\y  -\frac9{10}  s^{-\frac{17}{18}}e^{2\alpha s}\int_{B}|\grad w|^2 \rho
{\mathrm{d}}y\nonumber\\
&&+  s^{-\frac{17}{18}}e^{2\alpha s}\Sigma (s)  -\frac{17}{18}  s^{-\frac{35}{18}}e^{2\alpha s} {\cal M}(w(s))\nonumber\\
&&+(\frac{6}{5p+5}-\sigma_1)  s^{-\frac{17}{18}}e^{2\alpha s}\int_{B} |w|^{p+1}\frac{\rho}{1-|y|^2}\y\nonumber\\
&& +\frac{17-5p}{10p+10} s^{-\frac{17}{18}}e^{2\alpha s}\displaystyle\int_B 
|w|^{p+1}\rho{\mathrm{d}}y -\sigma_1s^{-\frac{17}{18}}e^{2\alpha s}\displaystyle\int_B 
w^{2}\rho{\mathrm{d}}y.
\end{eqnarray}
A simple integration over $[s, s']$ ensures that,   for  all $ s\ge \max(-\log T^*(x),1)$, we have 
\begin{eqnarray}\label{2411b}
{\cal F}_{1}(w(s'),s')\!\!&=&   \!\!
-\frac35 \int_{s}^{s'} \!\! {\tau}^{-\frac{17}{18}}e^{2\alpha {\tau}}\int_{B}|\nabla_{\theta}w|^2
\frac{|y|^2\rho}{1-|y|^2}\y
{\mathrm{d}}\tau  -\frac9{10}\int_{s}^{s'} \!\! {\tau}^{-\frac{17}{18}} e^{2\alpha {\tau}}\int_{B}|\grad w|^2 \rho
{\mathrm{d}}y
{\mathrm{d}}\tau\nonumber\\
&&+\int_{s}^{s'}  {\tau}^{-\frac{17}{18}}e^{2\alpha {\tau}}\Sigma ({\tau}) 
{\mathrm{d}}\tau -\frac{17}{18} \int_{s}^{s'}\! {\tau}^{-\frac{35}{18}}e^{2\alpha {\tau}} {\cal M}(w({\tau}))
{\mathrm{d}}\tau\nonumber\\
&&+(\frac{6}{5p+5}-\sigma_1)  \int_{s}^{s'} {\tau}^{-\frac{17}{18}}e^{2\alpha {\tau}}\int_{B} |w|^{p+1}\frac{\rho}{1-|y|^2}\y
{\mathrm{d}}\tau\nonumber\\
&& +\frac{17-5p}{10p+10} \int_{s}^{s'}{\tau}^{-\frac{17}{18}}e^{2\alpha {\tau}}\displaystyle\int_B 
|w|^{p+1}\rho{\mathrm{d}}y{\mathrm{d}}\tau\nonumber\\
&& -\sigma_1 \int_{s}^{s'}{\tau}^{-\frac{17}{18}}e^{2\alpha {\tau}}\displaystyle\int_B 
w^{2}\rho{\mathrm{d}}y
+{\cal F}_{1}(w(s),s).
\end{eqnarray}
To conclude the requested estimate, first we
control the third term to the right hand side.  
Clearly, by multiplying  \eqref{5dec123bis} by $s^{-\frac{17}{18}}e^{2\alpha s}$, we deduce after integrating over $[s, s']$ that
\begin{eqnarray}\label{2411c}
\int_{s}^{s'}  {\tau}^{-\frac{17}{18}}e^{2\alpha {\tau}}\Sigma ({\tau}) 
{\mathrm{d}}\tau  &\!\le &\!-\frac1{20}\int_{s}^{s'} \!{\tau}^{-\frac{17}{18}}e^{2\alpha {\tau}}\! \int_B \! (\partial_{s}w)^2\rho \y{\mathrm{d}}\tau+ \frac{4}{5} \int_{s}^{s'}\! {\tau}^{-\frac{17}{18}}e^{2\alpha {\tau}}\!\int_{B}|\grad w|^2\rho \y{\mathrm{d}}\tau  \nonumber\\
&&
+C_2\int_{s}^{s'} \!{\tau}^{-\frac{17}{18}}e^{2\alpha {\tau}}\int_{B} w^2\rho  \y{\mathrm{d}}\tau.
\end{eqnarray}
Also,  we need a  bound of the functional  ${\cal M}(w(s))$
  given by \eqref{Mdef}. From  the Hardy inequality  \eqref{Hardy},   we get
\begin{equation}\label{ml111}
 |{\cal M}(w(s))|
 \le C_3\displaystyle\int_B \Big((\partial_sw)^2+|\grad w|^2+w^2+
|w|^{p+1}\Big)\rho{\mathrm{d}}y.
\end{equation}
%
for some $C_3=C_3(N,p)>0$. Combining  \eqref{2411b},   \eqref{2411c}, and \eqref{ml111} we obtain
\begin{eqnarray}\label{2411d}
{\cal F}_{1}(w(s'),s')\!\!&\le&   \!  -(\frac1{10}-\frac{C_4}{s})\int_{s}^{s'} \!\! {\tau}^{-\frac{17}{18}} e^{2\alpha {\tau}}\int_{B}|\grad w|^2 \rho{\mathrm{d}}y{\mathrm{d}}\tau\nonumber\\
&& -(\frac1{20} -\frac{C_4}{s})\int_{s}^{s'} \!\! {\tau}^{-\frac{17}{18}} e^{2\alpha {\tau}}\int_{B}(\partial_s w)^2 \rho
{\mathrm{d}}y
{\mathrm{d}}\tau\nonumber\\
&&+(C_4-\sigma_1)  \int_{s}^{s'} {\tau}^{-\frac{17}{18}}e^{2\alpha {\tau}}\int_{B} |w|^{p+1}\frac{\rho}{1-|y|^2}\y
{\mathrm{d}}\tau\nonumber\\
&& +(C_4-\sigma_1) \int_{s}^{s'}{\tau}^{-\frac{17}{18}}e^{2\alpha {\tau}}\displaystyle\int_B 
w^2\rho{\mathrm{d}}y
{\mathrm{d}}\tau+{\cal F}_{1}(w(s),s),
\end{eqnarray}
for some $C_4=C_4(N,p)>0$.
By choosing $\sigma_1=C_4$,  and $S_1=S_1(N,p)>3$ large enough such that 
  \eqref{111} holds, 
 which ends the proof of  Proposition \ref{121}. \Box

\bigskip



\par

Thanks  to Proposition \ref{121} we get the following bounds:

\begin{cor}\label{coro33}
 For  all $s\ge \tilde {s_1}:= \max(2-\log T^*(x),S_1)$, we have  
\begin{equation}\label{K1}
0\le {\cal F}_{1}(w(s),s) 
\le {\cal F}_{1}(w(\tilde {s_1}),\tilde {s_1}),
\end{equation}
\begin{eqnarray}\label{K2}
\int_{\tilde  {s_1}}^{s}\tau^{-\frac{17}{18}}e^{2\alpha \tau}
 \int_{B}\Big((\partial_s w)^2 + |\grad w|^2\Big) \rho\y
{\mathrm{d}}\tau  &\le& 40 {\cal F}_{1}(w(\tilde {s_1}),\tilde {s_1}).
\end{eqnarray}
\end{cor}

\no {\it{Proof:}}

- {\it Proof of \eqref{K1}}: Note that the upper bound follows from  \eqref{111}. 
Therefore, it remains to prove the nonnegativity of the functional  ${\cal F}_{1}(w(s),s) $.
%
 By multiplying  the estimates  \eqref{ml111} by $s^{-\frac{17}{18}}e^{2\alpha s}$ in time between $s$ and $s+1$,  and using  (\ref{cor3}),  (\ref{F0}), and  (\ref{FF1}), we conclude that
\begin{eqnarray}\label{fin30}
\int_{s}^{s+1}  \tau^{-\frac{17}{18}}e^{2\alpha \tau}{\cal M}(w(\tau),\tau) 
d\tau\rightarrow 0,
 \ \ {\textrm {as}}\ \ s\rightarrow +\infty.
\end{eqnarray}
By exploiting  \eqref{LL}  and \eqref{fin30}, we deduce 
\begin{eqnarray}\label{fin300}
\int_{s}^{s+1}\!{\cal F}_{1}(w(\tau),\tau) 
d\tau\rightarrow 0,
 \ \ {\textrm {as}}\ \ s\rightarrow +\infty.
\end{eqnarray}
Combining the monotonicity of ${\cal  F}_{1}(w(s),s)$,
and  \eqref{fin300} we get 
\begin{eqnarray}\label{final1}
 {\cal F}_{1}(w(s),s) \rightarrow 0 ,\ \ \textrm{ as}\ \ \ s\rightarrow +\infty .
\end{eqnarray}
Finally, again by using  the monotonicity of ${\cal F}_{1}(w(s),s)$ and \eqref{final1}, we conclude
\begin{equation}\label{K1bis}
 {\cal F}_{1}(w(s),s) \ge0  \qquad \textrm{for  all } s\ge \tilde s_1=2+\max (-\log T^*(x),1).
\end{equation}
Then we deduce that  estimate (\ref{K1})
holds.

- {\it Proof of \eqref{K2}}:
 The estimates \eqref{K2} is  a direct
consequence of Proposition \ref{121} and \eqref{K1}. 
 This concludes the proof of Corollary \ref{coro33}.

\Box

\noindent{\it Proof of  Proposition \ref{improv} }
Clearly, by using (\ref{K2}) and the covering argument of \cite{MZimrn05},
we deduce the estimates \eqref{imp1}.
This concludes the proof of  Proposition \ref{improv}.
\Box


\bigskip

Now, we shall generalize the result obtained in the previous subsection under some additional hypotheses.  
Let us mention this result  is useful once we use the induction argument. 


\subsection{Slightly more
general improvments about  the  norm of $w$   in the space $\dot{H}^{1}(\R^N\times [-\log(s_0),\infty))$ under assumptions}\label{32}

In this subsection, we assume that there exist $k\in \N$ such that, we have the following:  

\begin{equation*}\label{HFF1}
(A_k)\qquad 
s^{\frac{k}{18}}e^{2\alpha s}\int_{s}^{s+1} \!\int_{B}\!\big((\partial_s
w(y,\tau))^2
+|\grad w(y,\tau)|^2\big)
{\mathrm{d}}y{\mathrm{d}}\tau\ \ \ \rightarrow 0
\textrm{ as}\ \ \ s\rightarrow +\infty,
\end{equation*}
\begin{equation*}\label{HF0}
(B_k)\qquad  \qquad \qquad \qquad  \qquad
\qquad \int_{\tilde s_0}^{\infty}s^{\frac{k}{18}}e^{2\alpha s}
\!\!\int_{B}\!|w(y,s)|^{p+1}{\mathrm{d}}y{\mathrm{d}}s \le
 \tilde K_k, \qquad\qquad
\end{equation*}
  where the constant $\tilde K_k$ depends only on $N, p, \delta_{0}(x_{0})$, 
 $T(x_{0}), k$ and 
$\|(u_0,u_1)\|_{H^{1}\times
L^{2}(B(x_0,\frac{T(x_0)}{\delta_0(x_0)}) )}$.
Furthermore,   we have
\begin{equation*}\label{HnesF0}
(C_k)\qquad \qquad \qquad \qquad \qquad \qquad 
s^{\frac{k}{18}}F_0(w(s),s)
\ \ \ \rightarrow 0
\textrm{ as}\ \ \ s\rightarrow +\infty.\qquad \qquad 
\end{equation*}

\begin{rem}
Let us recall the hypothesis the hypothesis $(A_0), (B_0)$ and  $(C_0)$ hold thanks to  \eqref{FF1}, \eqref{F0} and \eqref{nesF0}, where 
$\tilde K_0=K_2$.  
\end{rem}
Under the three  hypothesis,   $(A_k), (B_k)$ and  $(C_k)$, we are in position to prove the following:
\begin{pro}\label{Himprov}
Under the hypothesis $(A_k), (B_k)$ and  $(C_k)$, there exists $\tilde K_{k+1}>0,$ such that
 we have 
\begin{equation}\label{Himp1}
 \int_{\tilde s_0}^{\infty}\tau^{\frac{k}{18}-\frac{17}{18}}e^{2\alpha \tau}
   \int_{B}\Big( (\partial_s w)^2+ |\grad w|^2\Big) {\mathrm{d}}y
 {\mathrm{d}}\tau   \le \tilde K_{k+1}.
 \end{equation}  
\end{pro}

The proof follows a similar strategy to the proof of Proposition \ref{improv}. In fact,
in order to show Proposition \ref{Himprov}, we first
combine \eqref{cor3},$(A_k)$, $(B_k)$, and $(C_k)$, to  deduce the  following crucial estimate:
\begin{cor}\label{HP02244}
For all $s>\tilde s_0 $,  we have the following inequality
\begin{equation}\label{Hw2}
\int^{s}_{\tilde s_0}\tau^{\frac{k}{18}-\frac{17}{18}} e^{2\alpha \tau}
\int_{B}|w|^{p+1}\frac{\rho}{1-|y|^2} {\mathrm{d}}y{\mathrm{d}}s \leq  \Check K_{k+1}.
\end{equation}
\end{cor}

\noindent{\it Proof:} The proof is similar to  Corollary \ref{P02} by just using  $(A_k)$, $(B_k)$, and $(C_k)$, instead of 
 the estimates \eqref{cor3}, \eqref{F0}, and \eqref{FF1}.
 This concludes the proof of  Corollary \ref {HP02244}.
\Box

\par

Now, we shall  introduce the following functionals:
\begin{eqnarray}
{\cal U}_{k+1}(w(s),s)\!\!&=&\!\!\int_{s}^{\infty}\!\tau^{\frac{k}{18}-\frac{17}{18}}e^{2\alpha \tau}\int_{B}\! \Big(w^2\rho+|w|^{p+1}
\frac{ \rho}{1-|y|^2} \Big) {\mathrm{d}}y {\mathrm{d}}\tau,\qquad \label{HU}\\
{\cal F}_{k+1}(w(s),s)&= &s^{\frac{k}{18}-\frac{17}{18}}e^{2\alpha s} {\cal M}(w(s)) +\sigma_{k+1} {\cal U}_{k+1}(w(s),s),\label{HFF}
\end{eqnarray}
 for all $s\ge \tilde s_0,$ and where $\sigma_{k+1}$ is a positive constant chosen later.

\bigskip

Let us point out that 
the construction of a  family of   new Lyapunov
functional by induction argument in similarity variables is also an  important point of our
strategy which is also another novelty  in this paper.

 \begin{rem}
Let us point out   that, the  functional  ${\cal U}_{k+1}(w(s),s)$ given by \eqref{HU} is well defined.
This is a consequence of 
our previous result  given by \eqref{cor3}  and  the Corollary \ref{HP02244}. Furthermore, by using the estimates  \eqref{cor3} and \eqref{Hw2} we  can conclude that  
\begin{eqnarray}\label{HLL}
{\cal U}_{k+1 }(w(s),s) \rightarrow 0 ,\ \ \textrm{ as}\ \ \ s\rightarrow +\infty.
\end{eqnarray}
\end{rem}
Similarly, by following the strategy of  Proposition \ref{121}, 
 we can   choose $\sigma_{k+1}>0$,   such that  the function ${\cal F}_{k+1}(w(s),s)$
is  a decreasing functional for the equation
$\eqref{eqw1}$, for $s$ large enough. More precisely, we show the following proposition:
\begin{pro}\label{H12124}
There exists $S_{k+1}(N,p)>1$, such that 
for  all $s'\ge s\ge  \tilde s_{k+1}:= \max(2-\log T^*(x),S_{k+1})$, we have  the following 
\begin{equation}\label{H111}
{\cal F}_{k+1}(w(s'),s') +
\frac1{40} \int_{s}^{s'}\tau^{\frac{k}{18}-\frac{17}{18}}e^{2\alpha \tau}
   \int_{B}\Big( (\partial_s w)^2+ |\grad w|^2\Big) \rho {\mathrm{d}}y
 {\mathrm{d}}\tau   \le {\cal F}_{k+1}(w(s),s) .
 \end{equation}  
\end{pro}

%

\no {\it{Proof:}}
The argument is the same as  in Proposition \ref{121}.
  \Box
{ }

\par

\  \ \

Moreover,  thanks  to Proposition \ref{H12124} we get the following bounds:

\begin{cor}\label{Hcoro33}
 For  all $s\ge \tilde s_{k+1}:=\max(2-\log T^*(x),S_{k+1})$, we have  
\begin{equation}\label{HK1}
0\le {\cal F}_{k+1}(w(s),s) 
\le {\cal F}_{k+1}(w(\tilde s_{k+1}),\tilde s_{k+1}),
\end{equation}
\begin{eqnarray}\label{HK2}
\int_{{\tilde s_{k+1}}}^{s}\tau^{\frac{k}{18}-\frac{17}{18}}e^{2\alpha \tau}
 \int_{B}\Big((\partial_s w)^2 + |\grad w|^2\Big) \rho\y
{\mathrm{d}}\tau  &\le& 40 {\cal F}_{1}(w(\tilde s_{k+1}),\tilde s_{k+1}).
\end{eqnarray}
\end{cor}

\no {\it{Proof:}}
The argument is the same as in the corresponding part in Corollary \ref{coro33}.
  \Box

\noindent{\it Proof of  Proposition \ref{Himprov} }
Clearly, by using (\ref{HK2}) and the covering argument of \cite{MZimrn05},
we deduce the estimates \eqref{Himp1}.
This concludes the proof of  Proposition \ref{Himprov}.
\Box

\bigskip

\section{Proof of Theorem \ref{t0}}\label{sec4}
This section is devoted to the proof of
Theorem \ref{t0}. 
This section is divided into three parts:
\begin{itemize}
\item  In subsection 4.1,   based upon Proposition 1', we construct a  positive  Lyapunov functional 
for equation (\ref{C}).
\item In subsection 4.2, we prove   Theorem \ref{t0} in the particular case  $q=\frac1{18}$.
\item In subsection 4.3, we prove   Theorem \ref{t0} in the general case  by induction argument.
\end{itemize}
\subsection{Existence  of a Lyapunov functional for equation (\ref{C})}
 Consider $u $   a solution of (\ref{gen}) with
blow-up graph $\Gamma :\{x\mapsto T(x)\}$ and  $x_0$ is a non
characteristic point. 
We define the  following time:
$$t_1(x_0)=\max (T(x_0)-e^{-S_1},0) 
$$

Let   $T_0\in (t_1(x_0), T(x_0)]$, for all 
$x\in \er^N$  such that $|x-x_0|\le \frac{T_0}{\delta_0(x_0)}$, then we write $w$ instead of $w_{x,T^*(x)}$ defined in (\ref{scaling}) with $T^*(x)$ given in  (\ref{18dec1}).
Firstly,  for all
 $s\ge -\log (T^*(x))$  we recall the following natural 
functional:
\begin{eqnarray}
E_0(w(s))\!\!\!&=&\!\!\!\!\iint \Big(\frac{1}{2}(\partial_{s}w)^2+\frac{1}{2}(|\grad w|^2-(y\cdot \grad w)^2)+\frac{p+1}{(p-1)^2}w^2-\frac{|w|^{p+1}}{p+1}\Big) \y,\\\label{9sep1}
J_0(w(s))\!\!\!&=&\alpha \iint w\partial_{s}w \y-\frac{\alpha N}{2}
\!\!\iint w^2 \y,\\ 
F_0(w(s),s)\!\!\!&=&e^{2\alpha s}\big(E_0(w(s))+J_0(w(s))\big).
\end{eqnarray}

Moreover,  
   we define the functional
\begin{equation}\label{10dec2}
F_{1}(w(s),s)=s^{\frac1{18}}F_{0}(w(s),s)
+\frac1{18} \int_{s}^{\infty}\tau^{-\frac{17}{18}}F_{0}(w(\tau),\tau) {\mathrm{d}}\tau.
 \end{equation}

\begin{rem}
Let us   emphasize the fact that   the functional $ F_{1}(w(s),s)$ is well defined. In fact,  the estimates   \eqref{imp1} and  \eqref{cor3} implies 
\begin{equation*}
 \int_{\tilde s_0}^{\infty}\tau^{-\frac{17}{18}}e^{2\alpha \tau}
   \int_{B}(w^2+| \grad w|^2+
   (\partial_s w)^2){\mathrm{d}}y{\mathrm{d}}\tau
\le C.
 \end{equation*}  
By exploiting the basic inequality
\begin{equation}\label{b24}
 F_{0}(w(s),s)
 \le Ce^{2\alpha s}\displaystyle\int_B 
(w^2+| \grad w|^2+
   (\partial_s w)^2){\mathrm{d}}y,
\end{equation}
we deduce  the functional is defined and 
\begin{equation}\label{9sep10dec2}
\int_{s}^{\infty}\tau^{-\frac{17}{18}}F_{0}(w(\tau),\tau) {\mathrm{d}}\tau \rightarrow
0, \ \ {\textrm {as}}\  s\rightarrow +\infty.
 \end{equation}
\end{rem}

Now,
we derive that the functional  $F_{1}(w(s),s)$ is a decreasing 
  functional  of time  for equation (\ref{C}).
Let us first   state the following proposition:
\begin{pro}\label{l10dec}
\label{energylyap0} For all $s'\ge s\ge \max(-\log T^*(x),1)$, we have 
\begin{multline}\label{9sepAbis}
F_1(w(s'),s')+\int_{s}^{s'}\tau^{\frac{1}{18}}e^{2\alpha s} \int_{\partial B} \Big(\partial_s
w+\alpha  w\Big)^2 {\mathrm{d}}\sigma {\mathrm{d}}\tau\\ 
-\frac{\alpha (p-1)}{p+1}\int_{s}^{s'}s^{\frac1{18}} e^{2\alpha s}
\int_{B}|w|^{p+1}{\mathrm{d}}y{\mathrm{d}}\tau=F_1(w(s),s).
\end{multline}
Moreover, for all $s\geq max(
-\log (T^*(x),1)$, we have: 
\begin{equation}\label{26111}
F_{1}(w(s),s)\geq 0.
\end{equation}
\end{pro} 
{\it Proof}: 
Let us recall from Proposition 1 page 2320,   in \cite{HZdcds13}, 
the time derivative of the
functional  $F_0(w(s),s)$
is given by:
\begin{equation}\label{9sep00}
\frac{d}{ds}F_0(w(s),s)\;\;=\;\;-e^{2\alpha s} \int_{\partial B} \Big(\partial_s
w+\alpha  w\Big)^2 {\mathrm{d}}\sigma + \frac{\alpha (p-1)}{p+1}e^{2\alpha s}
\int_{B}|w|^{p+1}{\mathrm{d}}y.
\end{equation}
A straightforward computation yields the identity
\begin{equation}\label{9sep001}
\frac{d}{ds}F_1(w(s),s)=s^{\frac1{18}}\frac{d}{ds}F_0(w(s),s).
\end{equation}
From  \eqref{9sep00}, and \eqref{9sep001}  we get
\begin{multline}\label{9sepA}
\frac{d}{ds}F_1(w(s),s)=-s^{\frac1{18}}e^{2\alpha s} \int_{\partial B} \Big(\partial_s
w+\alpha  w\Big)^2 {\mathrm{d}}\sigma 
+ \frac{\alpha (p-1)}{p+1}s^{\frac1{18}}e^{2\alpha s}
\int_{B}|w|^{p+1}{\mathrm{d}}y.
\end{multline}
By integrating in time over $[s,s']$, the equality \eqref{9sepA}, we infer
 the desired estimate  \eqref{9sepAbis}. 
 Furthermore, by taking into
account the fact that $F_0(w(s),s)\ge0$, we conclude easily   
\eqref{26111},  which ends the proof of Proposition \ref{l10dec}.
\Box

\bigskip

Let us  use Proposition \ref{l10dec}. and the averaging technique of  and \cite{MZimrn05,MZma05} to get the following bounds:
\begin{lem}\label{lem1}
For  all $s\ge \hat s_0:= \max(-\log T^*(x),1)$, we have
\begin{equation}\label{F1}
0\le F_1(w(s),s)\le F_1(w(\hat s_0),\hat s_0),
\end{equation}
\begin{eqnarray}\label{F2}
\int_{\hat s_0}^{\infty}s^{\frac1{18}}e^{2\alpha s}
\!\!\int_{B}\!|w(y,s)|^{p+1}{\mathrm{d}}y{\mathrm{d}}s &\le&
-\frac{p+1}{\alpha (p-1)} F_1(w(\hat s_0),\hat s_0),
\end{eqnarray}
\begin{eqnarray}\label{F3}
\int_{\hat s_0}^{\infty}s^{\frac1{18}}e^{2\alpha s} \!\!\int_{\partial B}\!\Big(\partial
_sw(\sigma ,s)+\alpha w(\sigma ,s)\Big)^{2}{\mathrm{d}}\sigma
{\mathrm{d}}s &\le& F_1(w(\hat s_0),\hat s_0),
\end{eqnarray}
\begin{eqnarray}\label{corcor}
s^{\frac1{18}}e^{2\alpha s}\int_{s}^{s+1}\int_{B}\big (\partial_s w(y,\tau)\big)^2dy d\tau \rightarrow
0, \ \ {\textrm {as}}\  s\rightarrow +\infty.
\end{eqnarray}
\end{lem}
\no {\it{Proof:}}
  The first three estimates are  a direct
consequence of Proposition \ref{l10dec}. As for the last estimate, by introducing $f(y,s)=e^{\alpha s}w(y,s)$, we see that  the
dispersion estimate (\ref{F3}) can be written as follows:
\begin{eqnarray}\label{F4}
\int_{\hat s_0}^{\infty} \!\!\int_{\partial B}s^{\frac1{18}} \Big(\partial _sf(\sigma
,s)\Big)^{2}{\mathrm{d}}\sigma {\mathrm{d}}s &\le& F_1(w(\hat s_0),\hat s_0).
\end{eqnarray}
In particular, we have
\begin{eqnarray}\label{F5}
\int_{s}^{s+1} \!\tau^{\frac1{18}}\int_{\partial B}\!\Big(\partial _sf(\sigma
,\tau)\Big)^{2}{\mathrm{d}}\sigma {\mathrm{d}}\tau \rightarrow 0, \
\ {\textrm {as}}\ \ s\rightarrow +\infty.
\end{eqnarray}
By exploiting (\ref{F5}) where the space integration is done over the unit sphere, one can use the averaging technique of Proposition 4.2 (page 1147) in   \cite{MZimrn05} to get
\begin{eqnarray}\label{cor00001}
s^{\frac1{18}}e^{2\alpha s}\int_{s}^{s+1}\int_{B}\Big (\partial_s w(y,\tau
)-\lambda(\tau ,s) w(y,\tau )\Big )^2dy d\tau \rightarrow
0, \ \ {\textrm {as}}\  s\rightarrow +\infty,\ \
\end{eqnarray}
where $0\le \lambda(\tau ,s) \le C(\delta_0)$, for all $\tau \in
[s,s+1]$.
 By using  \eqref{cor3} and the above estimate, we conclude easily \eqref{corcor}
\Box


\medskip

\subsection{Proof of Theorem \ref{t0} in the particular case $q=\frac1{18}$}
\no Consider $u $ a solution of ({\ref{gen}}) with blow-up graph
$\Gamma:\{x\mapsto T(x)\}$. Translating Theorem \ref{t0} in the self-similar setting   
$w_{x_0,T_0}$ (we write $w$ for simplicity) defined
by (\ref{scaling}).  
%
First, 
from Lemma \ref{lem1}, we are in a position to  prove  the following:
\begin{pro}
\label{final}
If $u $   is a solution of ({\ref{gen}}) with blow-up graph
$\Gamma:\{x\mapsto T(x)\}$, then for all  $x_0\in \er^N$  and  $T_0\le
T(x_0)$, we have
\begin{eqnarray}\label{cor2}
s^{\frac1{18}}e^{2\alpha s}\int_{ B}\!|w(y,s)|^{\frac{p+3}2}{\mathrm{d}}y \ \
\rightarrow 0 ,\ \ \textrm{ as}\ \ \ s\rightarrow +\infty,
\end{eqnarray}
\begin{eqnarray}\label{cor3bis1}
s^{\frac{2 }{9(p+3)}}e^{\frac{8\alpha s}{p+3}}\int_{ B}\!|w(y,s)|^2{\mathrm{d}}y \ \
\rightarrow 0 ,\ \ \textrm{ as}\ \ \ s\rightarrow +\infty,
\end{eqnarray}
\begin{eqnarray}
s^{\frac1{18}}e^{2\alpha  s}\int_{s}^{s+1} \!\int_{B}\!
|\grad
w(y,\tau)|^2{\mathrm{d}}y{\mathrm{d}}\tau\ \ \ &\rightarrow& 0, \textrm{ as}\ \ \ s\rightarrow +\infty,\label{03}
\end{eqnarray}
Moreover, we have
\begin{eqnarray}\label{04}
s^{\frac1{18}}F_0(w(s),s)\ \rightarrow 0 ,\ \ \textrm{ as}\ \ \ s\rightarrow +\infty .
\end{eqnarray}
\end{pro}

\par

\no {\it{Proof 
:}}\\
\medskip

\no - {\it Proof of \eqref{cor2}}:
 Using the mean value theorem, we derive the existence of $\sigma(s)\in [s,s+1]$ such that
\begin{eqnarray}\label{s1}
s^{\frac1{18}}e^{2\alpha s}\int_{ B}\!|w(y,\sigma (s))|^{\frac{p+3}2}{\mathrm{d}}y
&=&\int_{s}^{s+1}\tau^{\frac1{18}}e^{2\alpha \tau}
\int_{ B}\!|w(y,\tau
)|^{\frac{p+3}2}{\mathrm{d}}y{\mathrm{d}}\tau.
\end{eqnarray}
By Jensen's inequality, we have
\begin{equation}\label{s2}
\int_{s}^{s+1}\tau^{\frac1{18}}e^{2\alpha \tau}\!\int_{ B}|w(y,\tau
)|^{\frac{p+3}2}{\mathrm{d}}y{\mathrm{d}}\tau \le  C \Big(\int_{s}^{s+1}\! \tau^{\frac1{18}}e^{2\alpha \tau}\int_{
B}\!\!|w(y,\tau
)|^{p+1}{\mathrm{d}}y{\mathrm{d}}\tau\Big)^{\frac{p+3}{2(p+1)}}.
\end{equation}
By combining  (\ref{s1}) and (\ref{s2}), we can write that
\begin{equation}\label{s3}
s^{\frac1{18}}e^{2\alpha s}\int_{ B}\!|w(y,\sigma (s))|^{\frac{p+3}2}{\mathrm{d}}y
\le C \Big(\int_{s}^{s+1} \tau^{\frac1{18}}e^{2\alpha \tau}\int_{ B}\!\!|w(y,\tau
)|^{p+1}{\mathrm{d}}y{\mathrm{d}}\tau\Big)^{\frac{p+3}{2(p+1)}}.
\end{equation}
Using  the fundamental theorem of calculus   together with (\ref{s3}), we deduce 
\begin{eqnarray*}
s^{\frac1{18}}e^{2\alpha s}\int_{ B}\!|w(y,s )|^{\frac{p+3}2}{\mathrm{d}}y&\le&
C\Big(\int_{s}^{s+1}\! \tau^{\frac1{18}}e^{2\alpha \tau}\int_{ B}\!\!|w(y,\tau
)|^{p+1}{\mathrm{d}}y{\mathrm{d}}\tau\Big)^{\frac{p+3}{2(p+1)}}\\
&& + C\Big(\int_{s}^{s+1}\tau^{\frac1{18}}e^{2\alpha \tau}\int_{
B}\!|w(y,\tau)|^{p+1}{\mathrm{d}}y{\mathrm{d}}\tau\Big)^{\frac12}\\
&& \quad \Big(\int_{s}^{s+1}\tau^{\frac1{18}}e^{2\alpha \tau}
\int_{ B}\!|\partial_s
w(y,\tau)|^2{\mathrm{d}}y{\mathrm{d}}\tau\Big)^{\frac12}.
\end{eqnarray*}
Since $\ds{s^{\frac1{18}} e^{2\alpha s}\int_{s}^{s+1}\int_{ B}\!|w(y,\tau
)|^{p+1}{\mathrm{d}}y{\mathrm{d}}\tau\rightarrow 0}$ from
(\ref{F2}), we use (\ref{corcor}) to obtain (\ref{cor2}).

\medskip

\no - {\it Proof of \eqref{cor3bis1}}: It follows from \eqref{cor2} through the Holder
inequality.

\medskip
%

\no - {\it Proof of \eqref{03}}:  We first prove  that, for all $\ge \hat s_0$,
\begin{eqnarray}\label{bb}
s^{\frac1{18}}e^{2\alpha s}\int_{s}^{s+1} \!\int_{B}|\grad
w(y,\tau)|^2{\mathrm{d}}y\le K.
\end{eqnarray}  
 By integrating  the functional $F_1(w,s)$  defined in  (\ref{10dec2}) in time between $s$ and $s+1$, we obtain:
\begin{eqnarray}\label{et}
&&\displaystyle\int_{s}^{s+1}\!\!\tau^{\frac1{18}}e^{2\alpha
\tau}\int_{B}\!\!\Big((\partial_sw)^2+|\grad w|^2-(y\cdot \grad w)^2\Big)
{\mathrm{d}}y{\mathrm{d}}\tau=  -\frac{2(p+1)}{(p-1)^2}
\displaystyle\int_{s}^{s+1}\!\! \tau^{\frac1{18}}e^{2\alpha \tau}\int_{B}\!\! w^2
{\mathrm{d}}y{\mathrm{d}}\tau\nonumber \\
&&
-\frac1{9}\int_{s}^{s+1} \int_{\tau}^{\infty}\zeta^{-\frac{17}{18}}F_1(w(\zeta ),\zeta )
d\zeta d\tau
+2\int_{s}^{s+1}\! F_1(w(\tau ),\tau )
d\tau\\
&&
+\frac{2}{p+1}\displaystyle\int_{s}^{s+1}\tau^{\frac1{18}}e^{2\alpha
\tau}\int_{B}\!\!|w|^{p+1} {\mathrm{d}}y{\mathrm{d}}\tau -\underbrace{2\alpha \int_{s}^{s+1}\! \tau^{\frac1{18}}e^{2\alpha
\tau}\displaystyle\int_{B}\!\!\big(w\partial_s w
 -\frac{
N}2 w^2\big) {\mathrm{d}}y{\mathrm{d}}\tau }_{A_0(s)}.\nonumber
\end{eqnarray}
Now, we    control all the terms on the right-hand side of the
relation (\ref{et}):

\no Note that the first two terms are negative, while the third term is
bounded because of the  bound (\ref{F1}) on the energy $F_1(w(s),s)$.
Remark that  (\ref{F2}) implies that the fourth term is also bounded.
Finally, it remains only to control the  term  $A_0(s)$.

\medskip

\no Combining  the Cauchy-Schwarz inequality, the inequality  $ab\le
 \nu a^2+ \frac1{4\nu}b^2$, and the fact that $N\ge 2$ and $\nu \in [0,1]$, we write
\begin{eqnarray}\label{et1}
 A_0(s)&\le& \frac{1}2 \int_{s}^{s+1}\!\!\tau^{\frac1{18}}e^{2\alpha
\tau}\displaystyle\int_{B}\!\!(\partial_s w)^2
{\mathrm{d}}y{\mathrm{d}}\tau.
\end{eqnarray}
By  combining  (\ref{et}), (\ref{et1}) and  the above-mentioned arguments for
the first four terms 
we get
\begin{eqnarray}\label{ett}
&&\displaystyle\int_{s}^{s+1}\!\!\tau^{\frac1{18}}e^{2\alpha
\tau}\int_{B}|\grad w|^2(1-|y|^2)
{\mathrm{d}}y{\mathrm{d}}\tau\le K.
\end{eqnarray}
By using (\ref{ett}) and the covering argument of \cite{MZimrn05},
we deduce the estimate \eqref{bb}.

Now, we are able to conclude  the improvement of the inequality \eqref{bb} by proving the estimate
(\ref{03}).

The first thing to  do is to recall from Lemma 2.2 page 2321,   in \cite{HZdcds13}, 
the time derivative of the
functional  $J_0(w(s))$
is given by:
\begin{eqnarray}\label{e08}
\frac{d}{ds}J_0(w(s))&=& \alpha \int_{B}(\partial_s w)^2{\mathrm{d}}y
- \alpha\int_{B}(|\grad w|^2-(y\cdot \grad w)^2){\mathrm{d}}y+\alpha \int_{B}|w|^{p+1}{\mathrm{d}}y\nonumber\\
&&-\alpha^2 \int_{\partial B}  w^2{\mathrm{d}}\sigma +(\alpha^2N  -\alpha \frac{2p+2}{(p-1)^2})\int_{B}w^2{\mathrm{d}}y\\
&&-2\alpha \int_{\partial B} w \partial_sw {\mathrm{d}}\sigma +2\alpha
\int_{B} (y\cdot\grad w)
\partial_sw {\mathrm{d}}y  -\alpha
(\frac{p+3}{p-1}-N) \int_{B}w\partial_sw{\mathrm{d}}y.\nonumber
\end{eqnarray}
Using \eqref{e08}, and remembering the definition of $J_0(w(s))$ and $\alpha$ given by \eqref{9sep1}
and \eqref{alpha} implies that
\begin{eqnarray}\label{e8}
\frac{d}{ds}J_0(w(s))&=& (-2\alpha- \frac1{18s} )J_0(w(s)) +\alpha \int_{B}(\partial_s w)^2{\mathrm{d}}y
-\alpha \int_{B}(|\grad w|^2-(y\cdot\grad w)^2){\mathrm{d}}y\nonumber\\
&&+\alpha \int_{B}|w|^{p+1}{\mathrm{d}}y
-\alpha^2 \int_{\partial B}  w^2{\mathrm{d}}\sigma -\alpha(
\frac{ N}{36s} 
 +\frac{2p+2}{(p-1)^2})\int_{B}w^2{\mathrm{d}}y\nonumber\\
&&-2\alpha \int_{\partial B} w \partial_sw {\mathrm{d}}\sigma +2\alpha
\int_{B} (y\cdot\grad w)
\partial_sw {\mathrm{d}}y  +\frac{\alpha }{18s}  \int_{B}w\partial_sw{\mathrm{d}}y.
\end{eqnarray}
Let $s\ge s_0+1$, $s_7=s_7(s)\in [s-1,s]$  and $s_8=s_8(s)\in
[s,s+1]$ to be chosen later. 
By integrating after multiplication by $s^{\frac1{18}}e^{2\alpha s}$ the
expression   \eqref{e8} in time between $s_7$ and $s_8$, we
obtain
\begin{eqnarray}\label{q2}
 &&-\alpha\int_{s_7}^{s_8}s^{\frac1{18}}e^{2\alpha s}\int_{B}|\grad w|^2(1-|y|^2){\mathrm{d}}y{\mathrm{d}}s=
\underbrace{s_8^{\frac1{18}}e^{2\alpha s_8}J_0(w(s_8))-s_7^{\frac1{18}}e^{2\alpha s_7}J_0(w(s_7))}_{A_1(s)}\qquad \qquad \nonumber\\
&&-\underbrace{2\alpha \int_{s_7}^{s_8}s^{\frac1{18}}e^{2\alpha s}\int_{B} (y\cdot\grad w)
\partial_sw {\mathrm{d}}y{\mathrm{d}}s}_{A_2(s)}
 -\underbrace{\alpha  \int_{s_7}^{s_8}s^{\frac1{18}}e^{2\alpha s}\int_{B}(\partial_s
w)^2{\mathrm{d}}y{\mathrm{d}}s}_{A_3(s)}\\
&&- \underbrace{\alpha \int_{s_7}^{s_8}s^{\frac1{18}}e^{2\alpha
s}\int_{B}|w|^{p+1}{\mathrm{d}}y{\mathrm{d}}s}_{A_4(s)} +\underbrace{
\alpha (\frac{ N }{36s} + \frac{2p+2}{(p-1)^2})
\int_{s_7}^{s_8}s^{\frac1{18}}e^{2\alpha s}\int_{B}w^2{\mathrm{d}}y{\mathrm{d}}s}_{A_5(s)}\nonumber\\
&&\underbrace{- \frac{\alpha}{18s}\int_{s_7}^{s_8} s^{\frac1{18}}e^{2\alpha
s}\int_{B}w\partial_sw{\mathrm{d}}y{\mathrm{d}}s}_{A_6(s)}\underbrace{-
\int_{s_7}^{s_8}s^{\frac1{18}}e^{2\alpha s}\int_{\partial B}  (\partial_sw)^2
{\mathrm{d}}\sigma{\mathrm{d}}s}_{A_7(s)}\nonumber\\  &&+
\underbrace{\int_{s_7}^{s_8}s^{\frac1{18}}e^{2\alpha s}\int_{\partial B}
(\partial_sw+\alpha w)^2{\mathrm{d}}\sigma
{\mathrm{d}}s}_{A_8(s)}\underbrace{+\alpha\int_{s_7}^{s_8}s^{\frac1{18}}e^{2\alpha
s}\int_{B}(|y|^2|\grad_{\theta} w|^2{\mathrm{d}}y{\mathrm{d}}s}_{A_{9}(s)}.\nonumber
\end{eqnarray}
Now, we    control all the terms on the right-hand
 side of the
relation (\ref{q2}):

\no Note that, by (\ref{9sep1}) and using the Cauchy-Schwarz
inequality, we can write
\begin{equation}\label{I1}
s_8^{\frac1{18}}e^{2\alpha s_8}|J_0(w(s_8))|\le Cs_8^{\frac1{18}}e^{2\alpha s_8}
\displaystyle\int_{B}(\partial_s w(s_8))^2 {\mathrm{d}}y+Cs_8^{\frac1{18}}e^{2\alpha
s_8} \displaystyle\int_{B}w^2(s_8)
 {\mathrm{d}}y.
\end{equation}
By exploiting  (\ref{cor3}) and the fact that  $s_8\in [s,s+1]$, we
conclude that
\begin{equation}\label{I11}
s_8^{\frac1{18}}e^{2\alpha s_8} \displaystyle\int_{B}w^2(s_8)
 {\mathrm{d}}y\rightarrow 0\ \ {\textrm {as}}\ \ s\rightarrow
 +\infty,
\end{equation}
on the one hand.  by   using the mean value
theorem, let us choose  $s_8=s_8(s)\in [s,s+1]$ such that
\begin{eqnarray}\label{q4}\displaystyle\int_{s}^{s+1}\tau^{\frac1{18}}
e^{2\alpha \tau}\displaystyle\int_{B}(\partial_s w(\tau))^2
{\mathrm{d}}y{\mathrm{d}}\tau&=&s_8^{\frac1{18}}e^{2\alpha
s_8}\displaystyle\int_{B}(\partial_s w(s_8))^2 {\mathrm{d}}y.
\end{eqnarray}
By combining  (\ref{corcor}) and (\ref{q4}) we
obtain
\begin{equation}\label{I111}
s_8^{\frac1{18}}e^{2\alpha s_8}\displaystyle\int_{B}(\partial_s w(s_8))^2
{\mathrm{d}}y\rightarrow 0\ \ {\textrm {as}}\ \ s\rightarrow
 +\infty.
\end{equation}
Then, by using (\ref{I1}), (\ref{I11}) and (\ref{I111}), we get
\begin{equation}\label{I1111}
s_8^{\frac1{18}}e^{2\alpha s_8}J_0(w(s_8))\rightarrow 0\ \ {\textrm {as}}\ \ s\rightarrow
+\infty.
\end{equation}
Similarly, we  can choose $s_7=s_7(s)\in [s-1,s]$ such that
\begin{equation}\label{I2}
s_7^{\frac1{18}}e^{2\alpha s_7}J_0(w(s_7))\rightarrow 0\ \ {\textrm {as}}\ \ s\rightarrow
+\infty.
\end{equation}
Note that by combining   (\ref{I1111}), and (\ref{I2}), 
we deduce that
\begin{equation}\label{B1}
A_1(s)\rightarrow 0\ \ {\textrm {as}}\ \ s\rightarrow +\infty.
\end{equation}
To estimate $A_2(s)$, since $s_7\in [s-1,s]$ and $s_8\in [s,s+1]$, we write
\begin{eqnarray}\label{q01}
|A_2(s)|&\le&C\Big(\int_{s-1}^{s+1}\tau^{\frac1{18}}e^{2\alpha \tau}\int_{B} |\grad w|^2
{\mathrm{d}}y{\mathrm{d}}\tau \Big)^{\frac12}
\Big(\int_{s-1}^{s+1}\tau^{\frac1{18}}e^{2\alpha \tau}\int_{B} (
\partial_sw)^2 {\mathrm{d}}y{\mathrm{d}}\tau \Big)^{\frac12}.\qquad
\end{eqnarray}
Then, by (\ref{q01}), (\ref{corcor}) and  (\ref{bb}), we deduce
\begin{equation}\label{B2}
A_2(s)\rightarrow 0\ \ {\textrm {as}}\ \ s\rightarrow +\infty.
\end{equation}
By (\ref{corcor}), we can say that
\begin{equation}\label{B3}
A_3(s)\rightarrow 0\ \ {\textrm {as}}\ \ s\rightarrow +\infty.
\end{equation}
By (\ref{F2}), we also deduce that
\begin{equation}\label{B4}
A_4(s)\rightarrow 0\ \ {\textrm {as}}\ \ s\rightarrow +\infty.
\end{equation}
Due to  (\ref{cor3}) and (\ref{corcor}), we have
\begin{equation}\label{B6}
A_6(s)\rightarrow 0\ \ {\textrm {as}}\ \ s\rightarrow +\infty.
\end{equation}
By (\ref{F3}), we write that
\begin{equation}\label{B8}
A_8(s)\rightarrow 0\ \ {\textrm {as}}\ \ s\rightarrow +\infty.
\end{equation}
Finally, 
by combining (\ref{B1}), (\ref{B2}), (\ref{B3}), (\ref{B4}), (\ref{B6}), (\ref{B8})
and the fact that the terms $A_5(s)$, $A_7(s)$ and $A_9(s)$  are negative, we conclude that
\begin{eqnarray}\label{q9}
 \int_{s}^{s+1}\tau^{\frac1{18}}e^{2\alpha \tau}\int_{B}|\grad w|^2(1-|y|^2){\mathrm{d}}y{\mathrm{d}}\tau\rightarrow 0
 \ \ {\textrm {as}}\ \ s\rightarrow +\infty.
\end{eqnarray}
By using (\ref{q9}) and the covering argument of \cite{MZimrn05},
we deduce that  estimate (\ref{03})
holds.


\medskip

\no - {\it Proof of \eqref{04}}:
By integrating the functional $F_1(w(s),s)$ defined in
(\ref{10dec2}) in time between $s$ and $s+1$, we  write
\begin{eqnarray}\label{fin2}
\int_{s}^{s+1}\!\!F_1(w(\tau),\tau)
d\tau&=&\displaystyle\int_{s}^{s+1}\int_{B}\!\!\tau^{\frac1{18}}e^{2\alpha  \tau}\Big (
\frac{1}{2}(\partial_s w)^2
+\frac{p+1}{(p-1)^2}w^2-\frac{1}{p+1}|w|^{p+1}\Big ) {\mathrm{d}}y{\mathrm{d}}\tau\nonumber\\
&&+\frac{1}{2}\displaystyle\int_{s}^{s+1}\! \tau^{\frac1{18}}e^{2\alpha 
\tau}\int_{B}\!\!\Big (|\grad w|^2-(y\cdot \grad w)^2\Big )
{\mathrm{d}}y{\mathrm{d}}\tau\\
&&+\alpha \int_{s}^{s+1}\!\!\tau^{\frac1{18}}e^{2\alpha 
\tau}\displaystyle\int_{B}w\partial_s w
{\mathrm{d}}y{\mathrm{d}}\tau- \frac{\alpha N}2
\int_{s}^{s+1}\!\tau^{\frac1{18}}e^{2\alpha  \tau}\displaystyle\int_{B}w^2
 {\mathrm{d}}y{\mathrm{d}}\tau\nonumber\\
 &&+\frac1{18}\int_{s}^{s+1} \int_{\tau}^{\infty}\zeta^{-\frac{17}{18}}F_{0}(w(\zeta),\zeta) {\mathrm{d}}\zeta{\mathrm{d}}\tau.\nonumber
\end{eqnarray}
By using (\ref{cor3}), \eqref{9sep10dec2},  (\ref{F2}), \eqref{corcor}, (\ref{03}), and  (\ref{fin2}), we conclude that
\begin{eqnarray}\label{fin3}
\int_{s}^{s+1}\!\!F_1(w,\tau)
d\tau\rightarrow 0,
 \ \ {\textrm {as}}\ \ s\rightarrow +\infty.
\end{eqnarray}
Combining the monotonicity of $F_1(w(s),s)$ proved in Proposition \ref{l10dec}, (\ref{fin3}), and \eqref{9sep10dec2}, we deduce  (\ref{04}).  This concludes
 the proof
 of  Proposition \ref{final}. Since Theorem \ref{t0} directly follows from Proposition \ref{final} through the self-similar change of variables \eqref{scaling}, this concludes the proof of Theorem \ref{t0} too in the particular case $q=\frac1{18}$.

\Box

\subsection{Proof of Theorem \ref{t0}}


Consider $u $   a solution of (\ref{gen}) with
blow-up graph $\Gamma :\{x\mapsto T(x)\}$ and  $x_0$ is a non
characteristic point. 
Let   $T_0\in (0, T(x_0)]$, for all 
$x\in \er^N$  such that $|x-x_0|\le \frac{T_0}{\delta_0(x_0)}$, then we write $w$ instead of $w_{x,T^*(x)}$ defined in (\ref{scaling}) with $T^*(x)$ given in  (\ref{18dec1}). \\

In this subsection, we state and proof the folowing:

\begin{pro}\label{lem124b}
For all $k\in \N$, for  all $s\ge \hat s_0:= \max(-\log T^*(x),1)$, we have
\begin{eqnarray}\label{F224b}
\int_{\hat s_0}^{\infty}s^{\frac{k}{18}}
e^{2\alpha s}
\!\!\int_{B}\!|w(y,s)|^{p+1}{\mathrm{d}}y{\mathrm{d}}s &\le&
\hat{K},
\end{eqnarray}
\begin{eqnarray}\label{corcor24b}
s^{\frac{k}{18}}e^{2\alpha s}\int_{s}^{s+1}\int_{B}\Big(\grad
w(y,\tau)|^2+\big (\partial_s w(y,\tau)\big)^2\big)dy d\tau \rightarrow
0, \ \ {\textrm {as}}\  s\rightarrow +\infty,
\end{eqnarray}
\begin{eqnarray}\label{0424b}
s^{\frac{k}{18}}F_0(w(s),s)\ \rightarrow 0 ,\ \ \textrm{ as}\ \ \ s\rightarrow +\infty .
\end{eqnarray}
\end{pro}

As stated in the introduction, the proof of  Theorem \ref{t0} is accomplished by induction argument. 
Let us mention that  
the beginning point of the induction argument the 
 hypothesis  $(A_0), (B_0)$ and  $(C_0)$ hold thanks to  \eqref{FF1}, \eqref{F0} and \eqref{nesF0}, where 
$\tilde K_0=K_2$ wich implies  the beginning point  $(k=0)$  of the induction argument is verified.
Then, we suppose that the three hypotheses,   $(A_k), (B_k)$ and  $(C_k)$ hold, for some $k\in \N$, and show that 
$(A_{k+1}), (B_{k+1})$ and  $(C_{k+1})$ valid  also.

%

Let us mention that  the  the hypothesis  $(A_1), (B_1)$ and  $(C_1)$ hold thanks to \eqref{corcor}, \eqref{F2},  and Proposition \ref{final}.
Naturally, we will use the same method 
to provide the second step of the induction argument.

Now, we define the functional:
\begin{equation}\label{10dec224}
F_{k+1}(w(s),s)=s^{\frac{k+1}{18}}F_{0}(w(s),s)
+\frac{k+1}{18} \int_{s}^{\infty}\tau^{\frac{k-17}{18}}F_{0}(w(\tau),\tau) {\mathrm{d}}\tau.\\
 \end{equation}

\begin{rem}
Let us   emphasize the fact that   the functional $ F_{k+1}(w(s),s)$ is well-defined. Indeed,  under the hypothesis $(A_k), (B_k)$ and  $(C_k)$,  the estimates   \eqref{Himp1} given in Proposition\ref{Himprov}, 
 \eqref{cor3}, and \eqref{b24},
we deduce  the functional  $F_{k+1}(w(s),s)$  is defined and 
\begin{equation}\label{9sep10dec224}
\int_{s}^{\infty}\tau^{\frac{k-17}{18}}F_{0}(w(\tau),\tau) {\mathrm{d}}\tau \rightarrow
0, \ \ {\textrm {as}}\  s\rightarrow +\infty.
 \end{equation}
\end{rem}

Now,
we derive that the functional  $F_{k+1}(w(s),s)$ is a decreasing 
  functional  of time  for equation (\ref{C}). More precisely, we easily show  the following:
\begin{pro}\label{l10dec24}
\label{energylyap0} Under the three hypotheses,   $(A_k), (B_k)$ and  $(C_k)$ , for all $s'\ge s\ge \max(-\log T^*(x),1)$, we have 
\begin{multline}\label{9sepAbis24}
F_{k+1}(w(s'),s')+\int_{s}^{s'}\tau^{\frac{k+1}{18}}e^{2\alpha s} \int_{\partial B} \Big(\partial_s
w+\alpha  w\Big)^2 {\mathrm{d}}\sigma {\mathrm{d}}\tau\\ 
-\frac{\alpha (p-1)}{p+1}\int_{s}^{s'}s^{\frac{k+1}{18}} e^{2\alpha s}
\int_{B}|w|^{p+1}{\mathrm{d}}y{\mathrm{d}}\tau=F_{k+1}(w(s),s).
\end{multline}
\end{pro} 
{\it Proof}: 
The proof follows a similar strategy to the proof of Proposition \ref{l10dec}.
\Box

\bigskip

Thanks to  Proposition \ref{l10dec24}. and the averaging technique of  and \cite{MZimrn05,MZma05} we get the following bounds:
\begin{lem}\label{lem124}
Under the three hypotheses,   $(A_k), (B_k)$ and  $(C_k)$, for  all $s\ge \hat s_0:= \max(-\log T^*(x),1)$, we have
\begin{equation}\label{F124}
0\le F_{k+1}(w(s),s)\le F_{k+1}(w(\hat s_0),\hat s_0),
\end{equation}
\begin{eqnarray}\label{F224}
\int_{\hat s_0}^{\infty}s^{\frac{k+1}{18}}
e^{2\alpha s}
\!\!\int_{B}\!|w(y,s)|^{p+1}{\mathrm{d}}y{\mathrm{d}}s &\le&
-\frac{p+1}{\alpha (p-1)} F_{k+1}(w(\hat s_0),\hat s_0),
\end{eqnarray}
\begin{eqnarray}\label{F324}
\int_{\hat s_0}^{\infty}s^{\frac{k+1}{18}}e^{2\alpha s} \!\!\int_{\partial B}\!\Big(\partial
_sw(\sigma ,s)+\alpha w(\sigma ,s)\Big)^{2}{\mathrm{d}}\sigma
{\mathrm{d}}s &\le& F_{k+1}(w(\hat s_0),\hat s_0),
\end{eqnarray}
\begin{eqnarray}\label{corcor24}
s^{\frac{k+1}{18}}e^{2\alpha s}\int_{s}^{s+1}\int_{B}\big (\partial_s w(y,\tau)\big)^2dy d\tau \rightarrow
0, \ \ {\textrm {as}}\  s\rightarrow +\infty.
\end{eqnarray}
\end{lem}
\no {\it{Proof:}}
  The first three estimates are  a direct
consequence of Proposition \ref{l10dec24}. Concerning  the last estimate,
 the proof follows a similar strategy to the similar part in  proof of Lemma \ref{lem1}.
\Box

\medskip

Thanks to Lemma \ref{lem124}, we are in a position to  prove  the following:
\begin{pro}
\label{final24}
If $u $   is a solution of ({\ref{gen}}) with blow-up graph
$\Gamma:\{x\mapsto T(x)\}$, then for all  $x_0\in \er^N$  and  $T_0\le
T(x_0)$, under the three hypotheses,   $(A_k), (B_k)$ and  $(C_k)$, we have
\begin{eqnarray}\label{cor224}
s^{\frac{k+1}{18}}e^{2\alpha s}\int_{ B}\!|w(y,s)|^{\frac{p+3}2}{\mathrm{d}}y \ \
\rightarrow 0 ,\ \ \textrm{ as}\ \ \ s\rightarrow +\infty,
\end{eqnarray}
\begin{eqnarray}\label{cor3bis124}
s^{\frac{2k+2 }{9(p+3)}}e^{\frac{8\alpha s}{p+3}}\int_{ B}\!|w(y,s)|^2{\mathrm{d}}y \ \
\rightarrow 0 ,\ \ \textrm{ as}\ \ \ s\rightarrow +\infty,
\end{eqnarray}
\begin{eqnarray}
s^{\frac{k+1}{18}}e^{2\alpha  s}\int_{s}^{s+1} \!\int_{B}\!
|\grad
w(y,\tau)|^2{\mathrm{d}}y{\mathrm{d}}\tau\ \ \ &\rightarrow& 0, \textrm{ as}\ \ \ s\rightarrow +\infty,\label{0324}
\end{eqnarray}
\begin{eqnarray}\label{0424}
s^{\frac{k+1}{18}}F_0(w(s),s)\ \rightarrow 0 ,\ \ \textrm{ as}\ \ \ s\rightarrow +\infty .
\end{eqnarray}
\end{pro}

\no {\it{Proof:}}
  The proof follows a similar method to the similar part in the proof of Proposition \ref{final}.
\Box

\no {\it{Proof of Proposition \ref{lem124b}:}} 
Clearly, by adding the fact that 
the beginning point of the induction argument the 
 hypothesis  $(A_0), (B_0)$ and  $(C_0)$ hold thanks to  \eqref{FF1}, \eqref{F0} and \eqref{nesF0}, and  Lemma\ref{lem124} and Proposition \ref{final24},
the induction is completed, 
 this concludes the proof of Proposition \ref{lem124b}.
\Box

\no {\it{Proof of  Theorem \ref{t0}:}} 
 Since Theorem \ref{t0} directly follows from Proposition  \ref{lem124b} through the self-similar change of variables \eqref{scaling}, this concludes the proof of Theorem \ref{t0}.

\Box

 \section{ Some elementary 
  identities  }
  In this section we prove two identities   which play an instrumental role
in our main argument.
First evaluate  the term
related to the Pohozaev multiplier given by:
\begin{equation}\label{A0}
{\cal{A}_{\e}}(s)=\int_{B}y\cdot\grad  w \Big(\div(\p \nabla w-\p (y\cdot\nabla w)y)\Big) \y,
\end{equation}
wehere $\e>0$.  More precisely, we prove   
the following identity:
\begin{lem}\label{L1} 
For all $w\in {\cal  H}$ it holds that
\begin{eqnarray}\label{mport01}
{\cal{A}_{\e}}(s)  &=&-\varepsilon \int_{B}|\nabla_{\theta}w|^2
\frac{|y|^2\p}{1-|y|^2}\y
-\varepsilon \int_{B}(y\cdot\grad  w)^2
\p{\mathrm{d}}y\noindent\\
&&+\frac{N}2\int_{B}\Big(|\grad w|^2-(y\cdot\grad  w)^2\Big) \p
{\mathrm{d}}y -\int_{B}|\grad w|^2 \p
{\mathrm{d}}y\cdot\nonumber
\end{eqnarray}
\end{lem}
{\it Proof}: 
We divide ${\cal{A}_{\e}}(s)  $ into three terms:
${\cal{A}_{\e}}(s)  ={\cal{B}_{\e}}(s)  +{\cal{C}_{\e}}(s)  +{\cal{D}_{\e}}(s),$ where
\begin{equation*}\label{A1}
{\cal{B}}_{\e}(s)=\int_{B}(y\cdot\grad  w) \Delta w\p \y,
\end{equation*}
\begin{equation*}\label{A2}
{\cal{C}}_{\e}(s)=-\int_{B}(y\cdot \grad  w) \div( (y\cdot \nabla w)y)\p \y,
\end{equation*}
and
\begin{equation*}\label{A0}
{\cal{D}_{\e}}(s)  =\int_{B}(y\cdot\grad  w) \Big( \nabla w-(y\cdot \nabla w)y \Big)\cdot \nabla \p \y .
\end{equation*}

To estimate ${\cal{A}}_{\e}(s)$, we start
observe the immediate  identity
\begin{equation}\label{A3}
  (y\cdot\grad w) \Delta w= \sum_{i,j}y_i\partial_i w
\partial^2_jw.
\end{equation}
By  integrating by parts, exploiting \eqref{A3} and  the fact that  $\ds{\sum_{i,j}\delta_{i,j}\partial_i w
\partial_jw=|\grad w|^2}$,      we can write
\begin{eqnarray}\label{A4}
{\cal{B}}_{\e}(s) &=&-\frac12\sum_{i,j}\int_{B}y_i\partial_i((
\partial_j w)^2)\p \y-\int_{B}|\grad w|^2
\p\y \no\\ &&-\sum_{i,j}\int_{B}y_i\partial_i w
\partial_j w\partial_j\p \y.\nonumber
\end{eqnarray}
By using the identity  $\partial_j\p=- \frac{2\varepsilon y_j}{1-|y|^2}\p$, 
\eqref{A4} and  integrating by part one has that
\begin{eqnarray}\label{A5}
{\cal{B}}_{\e}(s)  =\frac12\int_{B}|\grad w|^2\div (\p y)
{\mathrm{d}}y-\int_{B}|\grad w|^2 \p
{\mathrm{d}}y+2\varepsilon \int_{B}(y\cdot \grad w)^2\frac{\p}{1-|y|^2}
 \y.
\end{eqnarray}
Then,  from \eqref{A5}, and by using the identity  $\div (\p y)=N\p+y\cdot \nabla \p$,   we get 
\begin{eqnarray}\label{A6}
{\cal{B}}_{\e}(s)  =-\varepsilon \int_{B}|\grad w|^2\frac{|y|^2\p}{1-|y|^2}
{\mathrm{d}}y+2\varepsilon \int_{B}(y\cdot \grad w)^2\frac{\p}{1-|y|^2}\y+\frac{N-2}2\int_{B}|\grad w|^2 \p
{\mathrm{d}}y .\no
\end{eqnarray}
To estimate ${\cal{C}}_{\e}(s)$, we start with the use of the classical identity
\begin{equation}\label{A11}
\div( (y\cdot\nabla w)y)= N (y\cdot\nabla w)+ \grad (y\cdot\nabla w)\cdot y,
\end{equation}
and  integrating by part, to obtain
\begin{equation}\label{A12}
{\cal{C}}_{\e}(s) =-N\int_{B}(y\cdot\grad
w)^2\p{\mathrm{d}}y +\frac12 \int_{B}(y\cdot \grad w)^2 \div (
\p y) {\mathrm{d}}y.
\end{equation}
Similarly, by using a basic identity, 
we write
\begin{eqnarray}\label{A13}
{\cal{C}}_{\e}(s)&=&-\frac{N}2\int_{B}(y\cdot \grad  w)^2\p
{\mathrm{d}}y-\varepsilon \int_{B}(y\cdot \grad  w)^2
\frac{|y|^2\p}{1-|y|^2} {\mathrm{d}}y.
\end{eqnarray}
Furthermore, by exploiting  the identity $\grad \p=-\frac{2\e}{1-|y|^2}y,$  we conclude easily
\begin{equation}\label{A14}
{\cal{D}}_{\e}(s)
=
-2\e\int_{B}(y\cdot\grad  w)^2 \p\y.
\end{equation}
By combining  \eqref{A6}, \eqref{A13} and  \eqref{A14} 
we deduce easily
\begin{eqnarray}
{\cal{A}}(s)  &=&-\varepsilon \int_{B}(|\grad w|^2-(y\cdot\grad  w)^2)
\frac{|y|^2\p}{1-|y|^2}
{\mathrm{d}}y\\
&&+\frac{N}2\int_{B}(|\grad w|^2-(y\cdot\grad  w)^2)\p
{\mathrm{d}}y-\int_{B}|\grad w|^2 \p
{\mathrm{d}}y.\nonumber
\end{eqnarray}
From the identity $|\grad w|^2-(y\cdot\grad  w)^2)=|\grad_{\theta} w|^2+(1-|y|^2)|\grad_r  w|^2$, the above estimate implies
\eqref{mport1}, which ends the proof of Lemma \ref{L1}. 
\Box

\bigskip

Then, we evaluate  the term
\begin{equation}\label{A9}
{\cal{E}_{\e}}(s)=- \int_{B}\div(\p \grad w-\p (y\cdot\grad w)
y)w\frac1{\sqrt{1-|y|^2}}{\mathrm{d}}y,
\end{equation}
  More precisely, we prove   
the following identity:
\begin{lem}\label{L52} 
For all $w\in {\cal  H}$ it holds that
\begin{eqnarray}\label{mport1}
{\cal{E}_{\e}}(s)
 &=& \int_{B}|\grad w_{\theta}|^2\pp{\mathrm{d}}y+\int_{B}|\grad w_r|^2{\rho_{\e+\frac12}}{\mathrm{d}}y
\\
&&+\int_{B}w y\cdot\grad w
 \frac{\p}{\sqrt{1-|y|^2}}{\mathrm{d}}y.\nonumber
\end{eqnarray}
\end{lem}
{\it Proof}: 
By integrating by parts, and using the fact that
$$ \int_{\partial B}y\cdot (\p \grad w-\p (y\cdot\grad w)
y)w\frac1{\sqrt{1-|y|^2}}{\mathrm{d}}\sigma= \int_{\partial B}(y\cdot \grad w) 
w{\sqrt{1-|y|^2}}\p {\mathrm{d}}\sigma=0,$$
 we write
\begin{equation}\label{AA1}
{\cal{E}_{\e}}(s)
 =\int_{B}(\p \grad w-\p (y\cdot \grad w)
y)\nabla (w\frac1{\sqrt{1-|y|^2}}){\mathrm{d}}y.
\end{equation}
A straightforward computation yields
 the identity
\begin{equation}\label{AA0}
\nabla (w\frac1{\sqrt{1-|y|^2}})=\frac1{\sqrt{1-|y|^2}} \nabla w+ \frac{1}{(1-|y|^2)^{\frac32}}wy.
\end{equation}
Substituting \eqref{AA0} into \eqref{AA1}, we have that
\begin{eqnarray}\label{V1}
{\cal{E}_{\e}}(s)
 &=&\int_{B}\Big(|\grad w|^2-(y\cdot\grad w)^2\Big)\frac{\p}{\sqrt{1-|y|^2}}{\mathrm{d}}y\nonumber\\
&&+\int_{B}w y\cdot\grad w
 \frac{\p}{\sqrt{1-|y|^2}}{\mathrm{d}}y.
\end{eqnarray}
By  using \eqref{wr2}, we write easily
\begin{equation}\label{V2}
\int_{B}\Big(|\grad w|^2-(y\cdot\grad w)^2\Big)\frac{\p}{\sqrt{1-|y|^2}}{\mathrm{d}}y
  =\int_{B}|\grad w_{\theta}|^2\pp{\mathrm{d}}y+\int_{B}|\grad w_r|^2{\rho_{\e+\frac12}}{\mathrm{d}}y.
\end{equation}
Substituting  \eqref{V2} into \eqref{V1}, we have 
\begin{eqnarray}\label{V3}
{\cal{E}_{\e}}(s)
 &=& \int_{B}|\grad w_{\theta}|^2\pp{\mathrm{d}}y+\int_{B}|\grad w_r|^2{\rho_{\e+\frac12}}{\mathrm{d}}y
\\
&&+\int_{B}w y\cdot\grad w
 \frac{\p}{\sqrt{1-|y|^2}}{\mathrm{d}}y.\nonumber
\end{eqnarray}
This ends the proof of Lemma \ref{L52}. 
\Box

\def\cprime{$'$}

\noindent{\bf Address}:\\
 Department of Basic Sciences, Deanship of Preparatory and Supporting Studies,
  Imam Abdulrahman Bin Faisal University
P.O. Box 1982 Dammam, Saudi Arabia.\\
\vspace{-7mm}
\begin{verbatim}
e-mail:  mahamza@iau.edu.sa
\end{verbatim}
Universit\'e Paris 13, Institut Galil\'ee,
Laboratoire Analyse, G\'eom\'etrie et Applications, CNRS UMR 7539,
99 avenue J.B. Cl\'ement, 93430 Villetaneuse, France.\\
\vspace{-7mm}
\begin{verbatim}
e-mail: Hatem.Zaag@univ-paris13.fr
\end{verbatim}

\end{document}